\documentclass[lettersize,journal]{IEEEtran}
\usepackage{amsmath,amsfonts}
\usepackage{algorithmic}
\usepackage{array}
\usepackage{textcomp}
\usepackage{stfloats}
\usepackage{url}
\usepackage{verbatim}
\usepackage{graphicx}

\usepackage{cellspace}

\usepackage{booktabs}
\usepackage{amssymb}
\usepackage{cite}
\usepackage{epsfig}
\usepackage{epstopdf}
\usepackage{color}
\usepackage{bm}
\usepackage{hyperref}

\usepackage{enumitem}
\usepackage{xspace}
\usepackage{microtype}
\usepackage{latexsym,amsthm,graphics,amscd,verbatim,mathrsfs}
\usepackage{algorithm}
\usepackage{algorithmic}

\newcommand{\norm}[1]{\|#1\|}
\usepackage{orcidlink}
\usepackage{caption}
\usepackage{subcaption}
\usepackage{tabularx, boldline}

\def\BibTeX{{\rm B\kern-.05em{\sc i\kern-.025em b}\kern-.08em T\kern-.1667em\lower.7ex\hbox{E}\kern-.125emX}}

\DeclareMathOperator{\dist}{dist}
\newcommand{\ie}{{\it i.e.}}
\newcommand{\alg}{VR$^3$PM\xspace}

\newtheorem{assumption}{Assumption}
\newtheorem{theorem}{Theorem}
\newtheorem{lemma}{Lemma}
\newtheorem{proposition}{Proposition}

\theoremstyle{definition}
\newtheorem{example}{Example}

\usepackage{balance}

\begin{document}
\title{Variance Reduced Random Relaxed Projection Method for Constrained Finite-sum  Minimization Problems}

\author{ 
	Zhichun Yang, Fu-quan Xia, Kai Tu$^{\orcidlink{0000-0002-0557-6792}}$, and Man-Chung Yue$^{\orcidlink{0000-0002-7992-9490}}$
\thanks{This work was supported in part by National Natural Science Foundation of China (No. 12101436 and No. 12271217) and the Hong Kong Research Grants Council under the GRF project 15305321.
Zhichun Yang (yangzhichun1994@163.com) and Fu-quan Xia (fuquanxia@163.com) are with the School of Mathematical Science, Sichuan Normal University, Chengdu, China. 
Kai Tu (kaitu\_02@163.com, corresponding author) is with School of Mathematical Sciences, Shenzhen University, Shenzhen, China. Man-Chung Yue (mcyue@hku.hk) is with Musketeers Foundation Institute of Data Science and Department of Industrial and Manufacturing Systems Engineering, The University of Hong Kong, Hong Kong, China.}
}

\markboth{}%
{}

\maketitle

\begin{abstract}
For many applications in signal processing and machine learning, we are tasked with minimizing a large sum of convex functions subject to a large number of convex constraints. In this paper, we devise a new random projection method (RPM) to efficiently solve this problem. Compared with existing RPMs, our proposed algorithm features two useful algorithmic ideas. First, at each iteration, instead of projecting onto the subset defined by one of the constraints, our algorithm only requires projecting onto a half-space approximation of the subset, which significantly reduces the computational cost as it admits a closed-form formula. Second, to exploit the structure that the objective is a sum, variance reduction is incorporated into our algorithm to further improve the performance. As theoretical contributions, under a novel error bound condition and other standard assumptions, we prove that the proposed RPM converges to an optimal solution and that both optimality and feasibility gaps vanish at a sublinear rate. In particular, via a new analysis framework, we show that our RPM attains a faster convergence rate in optimality gap than existing RPMs when the objective function has a Lipschitz continuous gradient, capitalizing the benefit of the variance reduction. We also provide sufficient conditions for the error bound condition to hold. Experiments on a beamforming problem and a robust classification problem are also presented to demonstrate the superiority of our RPM over existing ones.
\end{abstract}

\begin{IEEEkeywords}
Constrained Optimization, Finite-Sum Minimization, Random Projection Method, Relaxed Projection, Variance Reduction, Convergence Rate Analysis.
\end{IEEEkeywords}

\section{Introduction}\label{sec:intro}
\IEEEPARstart{T}{his} paper considers the following constrained convex optimization problem:
\begin{equation}
	\label{1.1}
	\begin{array}{c@{\quad}l}
		\displaystyle\min_{\bm{x}} & \displaystyle f(\bm{x})=\frac{1}{n}\sum_{i =1}^n f_i(\bm{x}) \\
		\noalign{\smallskip}
		\mbox{s.t.} & \bm{x}\in C = C_0\cap C_{[m]}, 
	\end{array}
\end{equation}
where $f_i:\mathbb{R}^d\rightarrow \mathbb{R}$ is a differentiable convex function for $i=1,\dots,n$, $C_0\subseteq \mathbb{R}^d$ is a non-empty, closed and convex set, $C_{[m]} = \cap_{j\in [m]} C_j$, and $C_j=\{\bm{x}\in \mathbb{R}^{d}\mid \phi_j(\bm{x})\leq 0\}$ with $\phi_j : \mathbb{R}^d\rightarrow \mathbb{R}$ being a convex but possibly non-differentiable function for $j \in [m]$. Here and throughout the paper, $[t] = \{1,\dots,t\}$ for any positive integer $t$.
Problem~\eqref{1.1} finds applications across a wide range of areas. Below we present two concrete examples.

\begin{example}\label{exam:constr_LASSO}
The first example is (linearly) constrained least squares problem:
\begin{equation}
\label{opt:constr_LASSO}
	\begin{array}{c@{\quad}l}
		\displaystyle\min_{\bm{x}} & \displaystyle \| A \bm{x} - \bm{a} \|^2  \\
		\noalign{\smallskip}
		\mbox{s.t.} & B \bm{x} \le \bm{b}, 
	\end{array}
\end{equation}
where $A\in \mathbb{R}^{\ell \times d}$, $\bm{a}\in \mathbb{R}^\ell$, $B\in\mathbb{R}^{m\times d}$ and  $\bm{b}\in \mathbb{R}^m$.
Problem~\eqref{opt:constr_LASSO} subsumes constrained LASSO~\cite{gaines2018algorithms, deng2020efficient} and generalized LASSO~\cite{xu2013generalized} as special cases and can be applied to many signal processing problems, including linear quadratic state estimation~\cite{kalman1960new}, trend filtering~\cite{kim2009ell_1}, wavelet smoothing~\cite{donoho1995adapting} and neural decoding~\cite{ng2010generalized}.
\end{example}

\begin{example}\label{exam:DRC}
Another example of problem~\eqref{1.1} arises from the context of robust classification and will be used in our numerical experiments in Section~\ref{sec:experiment}. To begin, we define a distance between two probability distributions $\mathbb{P}_1$ and $\mathbb{P}_2$:
\begin{align*}
	&W (\mathbb{P}_1, \mathbb{P}_2) \\
	= &\min_{\mathbb{Q}\in \Gamma(\mathbb{P}_1, \mathbb{P}_2)} \int_{ \Xi\times \Xi } \!\left( \| \bm{w}_1 - \bm{w}_2 \|_\infty + |y_1 - y_2|\right)  \mathrm{d}\mathbb{Q}((\bm{w}_1,y_1), (\bm{w}_2,y_2)),
\end{align*}
where $\Xi = \mathbb{R}^\ell \times \{+1, -1\}$ and $\Gamma(\mathbb{P}_1, \mathbb{P}_2)$ is the set of joint distributions on $\Xi\times \Xi$ with the first marginal $\mathbb{P}_1$ and second marginal $\mathbb{P}_2$. 
The distance $W$ is an instance of Wasserstein distance~\cite{yue2022linear}. 
The ball centered at $\widehat{\mathbb{P}}$ of radius $\epsilon > 0$ in the space of probability distributions is denoted by $B_\epsilon(\widehat{\mathbb{P}})$. Based on ideas from distributionally robust optimization~\cite{kuhn2019wasserstein}, the following model for robust binary classification is considered in~\cite{shafieezadeh2015distributionally}:
\begin{equation*}
	\label{opt:DRO}
	\begin{array}{c@{\quad}l}
		\displaystyle\min_{\bm{u}\in \mathbb{R}^\ell}\sup_{\mathbb{P} \in B_\epsilon(\widehat{\mathbb{P}})}  \mathbb{E}_{(\bm{w}, y)\sim \mathbb{P}}\left[ \log\left( 1 + \exp\left( - y\, \bm{u}^\top \bm{w} \right) \right) \right] ,
	\end{array}
\end{equation*}
which is intractable in general, due to the infinite-dimensionality of the maximization. However, if the center $\widehat{\mathbb{P}}$ is an empirical distribution associated with a sample $\{ (\bm{w}_i, y_i) \}_{i \in [n]}$, then by~\cite{shafieezadeh2015distributionally} it is equivalent to the problem
\begin{equation}
	\label{opt:DRLR}
	\begin{array}{c@{\quad}l}
		\displaystyle\min_{\bm{u},\bm{s},\lambda} & \displaystyle\lambda \epsilon +\frac{1}{n} \! \sum_{i \in [n]}\left( s_i+\log\left( 1 + \exp\left( - y_i \, \bm{u}^\top \bm{w}_i \right) \right)\right) \\
		\noalign{\smallskip}
		\mbox{s.t.} &  y_j\,\bm{u}^{\top}\bm{w}_j  \le s_j + \lambda, \  j \in [n]\\
		& \|\bm{u}\| \le \lambda,\ \bm{u}\in \mathbb{R}^\ell ,\, \bm{s}\in\mathbb{R}_+^n ,\, \lambda \in \mathbb{R}.
	\end{array}
\end{equation}
\end{example}
\noindent In the notation of problem~\eqref{1.1}, $d = \ell + 1 + n$, $m = n$, $\bm{x} = (\bm{u}^\top, \lambda, \bm{s}^\top)^\top$, $C_0 = \{ (\bm{u}^\top, \lambda, \bm{s}^\top)^\top: \|\bm{u}\| \le \lambda,\ \bm{s}\in\mathbb{R}_+^n \}$,  $f_i (\bm{x}) = \epsilon \lambda + s_i + \log\left(1 + \exp\left( -y_i\, \bm{u}^\top \bm{w}_i \right)\right)$ and $\phi_j (\bm{x}) = y_j\,\bm{u}^{\top}\bm{w}_j - \lambda  - s_j$ for $i,j\in[n]$.

A straightforward algorithm for problem~\eqref{1.1} is the projected gradient method~(PGM) whose iteration takes the form
$
\bm{x}^{k+1}=\Pi_C \left(\bm{x}^k-\alpha_k\nabla f(\bm{x}^k) \right)
$,
where $\alpha_k> 0$ is the step-size and $\Pi_C(\cdot)$ denotes the projection map onto $C$.
The theory on PGM is rather complete, at least for convex problems. For example, it can be proved that under mild assumptions, PGM converges to an optimal solution with a sublinear rate $\mathcal{O}(1/K)$~\cite{beck2017first}, where $K$ is the total number of iterations. Under stronger assumptions, it is also proved that PGM converges linearly to an optimal solution~\cite{Nesterov2018}  (\ie, $\mathcal{O}(c^K)$ for some $c\in (0,1)$). 

However, PGM is not a viable algorithm for solving problem~\eqref{1.1} if
\begin{enumerate}[label=(\roman*)]
	\item \label{difficulty_i} the number $n$ of summands $f_i$ is large, 
	
	\item \label{difficulty_ii} the number $m$ of constraints is large, or
	
	\item \label{difficulty_iii} the projections onto some of the subsets $C_j$ are difficult.
\end{enumerate}
In case of difficulty~\ref{difficulty_i}, it is expensive to compute the gradient $\nabla f(\bm{x}^k)$; and when we have difficulty~\ref{difficulty_ii} or~\ref{difficulty_iii}, computing the projection $\Pi_C$ onto the whole feasible region $C$ is highly computationally demanding, if not impossible.

A standard idea to handle difficulty~\ref{difficulty_i} is to replace the gradient $\nabla f(\bm{x}^k)$ by the random estimator $\nabla f_{i_k} (\bm{x}^k)$, where $i_k$ is a random index. The resulting algorithm is known as the stochastic gradient method (SGM), which traces back to the work~\cite{SQG27} 
of Robbins and Monro.
A natural generalization of SGM is the mini-batch SGM~\cite{bottou2018optimization}, where instead of a single summand $f_i$, multiple summands are used to form the random gradient estimator.
SGM and its variants have become arguably the most popular algorithms for modern, large-scale machine learning. 
One particular reason comes from its significantly lower per iteration cost, compared to that of PGM. For more details, we refer the readers to the recent survey~\cite{Net2019}.
A major drawback of SGM and its mini-batch variant is that its gradient estimator tends to introduce a large variance to the algorithm, which necessitates the use of conservative step-size and in turn leads to slow convergence. Indeed, the convergence rate of SGM for minimizing smooth~(non-strongly) convex function is only $\mathcal{O}(1/\sqrt{K})$~\cite{nemirovski2009robust}, which is worse than the rate $\mathcal{O}(1/K)$ of the standard PGM based on the exact gradient. Similarly, to minimize a smooth strongly convex function, SGM can only achieve a slower rate of $\mathcal{O}(1/K)$~\cite{bottou2018optimization}, as compared to the linear convergence rate of the gradient descent.
Various variance reduction techniques have been proposed to remedy the variance issue. Notable examples of variance reduced SGM include SAG~\cite{ roux2012stochastic}, SAGA~\cite{AVR15} and SVRG~\cite{AVS21}. Thanks to the variance reduction techniques, the convergence rates of these variants match with that of the gradient descent on both non-strongly and strongly convex problems.

To deal with optimization problems with a large number of constraints, \ie, in the presence of difficulty~\ref{difficulty_ii}, one could adopt the framework of random projection methods~(RPMs). Roughly speaking, at each iteration, an RPM improves feasibility with respect to the intersection~$C_{[m]}$ by using only one randomly selected subset $C_{j_k}$ but ignores the others, where $j_k$ is a random index. The computational cost at each iteration can thus be substantially reduced.
In the context of optimization, RPMs were first studied by Nedi\'c in~\cite{Nedic}, which can be viewed as an extension of the algorithm in~\cite{polyak2001random} by Polyak for convex feasibility problems to convex minimization problems. Subsequent works on RPMs include~\cite{wang2015random, wang2016stochastic, nedic2019, Necoara2022}. On the theoretical side, these works typically show that under certain assumptions, both the optimality and feasibility gaps converge to zero at a sublinear rate.

Some of these RPMs, such as~\cite{wang2015random, wang2016stochastic, Necoara2022}, are designed for constrained optimization where the objective is an expectation $f(\bm{x}) = \mathbb{E}_{I} [ f_I (\bm{x} ) ]$ and subsumes the finite-sum objective in problem~\eqref{1.1} as a special case. In these works, it is assumed that given any $\bm{x}$, one can access~$\nabla f_i (\bm{x})$ at some realization $i$ of the random variable~$I$, serving as a random gradient estimator, which is used in the updating step of the algorithms. Therefore, it similarly suffers from the variance issue. In view of our previous discussion, it is tempted to replace the gradient estimator~$\nabla f_i (\bm{x}^k)$ by its variance reduced counterpart. Unfortunately, it is unclear whether the variance reduced estimators satisfy the assumptions required by the existing RPMs, such as~\cite[Assumption~1(c)]{wang2016stochastic} and \cite[Assumption~1]{Necoara2022}. Hence, incorporating variance reduction into RPMs with provable theoretical guarantees is highly nontrivial.
One of the contributions of this paper is to show that an improved convergence rate can be obtained by integrating the SVRG variance reduction technique into RPMs.

When difficulty~\ref{difficulty_iii} is present, computing the projection onto the subset $C_j$ could be time consuming. This hinders the the practicality of projection-based RPMs, such as~\cite{wang2015random, wang2016stochastic}, where each iteration requires computing the projection $\Pi_{C_{j_k}}$ onto the randomly selected subset $C_{j_k}$.
In such a case, one could improve the efficiency by approximating the subset $C_{j_k}$ by a half-space at each iteration. The advantage is that the projection onto a half-space can be computed much more efficiently via an explicit formula. The idea of approximating complicated feasible region by half-spaces is not new. It has been studied in many other settings and under different names, including the outer approximation method~\cite{Fuku84} and the cutting-plane method~\cite{Kelley}.

In this paper, we develop a new RPM that aims at solving problem~\eqref{1.1} in the face of difficulties~\ref{difficulty_i}-\ref{difficulty_iii}.
The proposed algorithm features two useful algorithmic ideas that can significantly improve the practical performance: variance reduction and half-space approximation of the complicated subsets. To the best of our knowledge, this is the first time these two ideas are simultaneously incorporated into the framework of RPMs. Furthermore, the proposed RPM enjoys rigorous theoretical guarantees. Under assumptions similar to previous works, we prove that the sequence of iterates generated by our RPM converges to an optimal solution to problem~\eqref{1.1}. 
Moreover, the convergence rates of the optimality gap and the feasibility gap are $\mathcal{O}(1/\sqrt{K})$ and $\mathcal{O}(1/K)$, respectively. If $f$ has a Lipschitz gradient, then the advantage of the variance reduction can be realized theoretically, and the rate of the optimality gap is improved to $\mathcal{O}(1/K^{2/3})$. Finally, under a quadratic growth condition on $f$, the rates can be strengthened to $\mathcal{O}(1/K)$ and $\mathcal{O}(1/K^2)$, respectively.
Because of the integration of the two new algorithmic ideas, theoretical frameworks of existing RPMs do not apply, and new analyses are developed.
In particular, we have formulated and proved an error bound-type condition, which plays an instrumental role in the theoretical development of our RPM. We believe that it could be of independent interest in the study of other optimization algorithms.

The rest of this paper is organized as follows. 
We present the proposed algorithm and the required assumptions in Section~\ref{sec:VR3PM} and Section~\ref{sec:assumptions}, respectively.
In Section~\ref{subsec2}, we study the convergence behaviour of the proposed algorithm.
Numerical results are reported in Section~\ref{sec:experiment}.


\subsection{Notation}\label{sec:notation}
The Euclidean norm and inner product  are denoted by $\|\cdot\|$ and $\langle\cdot,\cdot \rangle$, respectively.
For two non-negative sequences $\{a_k\}$ and $\{b_k\}$, we write $a_k=\mathcal{O}(b_k)$ if there exists $c>0$ such that $a_k\leq c\,b_k$ for any $k\geq 0$.
We denote by $|S|$ the cardinality of a set $S$.
For any convex function $\phi$ and $\bm{x}\in \text{dom}(\phi)$, we denote by $\partial \phi(\bm{x}) = \{ \bm{\xi}\in \mathbb{R}^{d}\mid \phi(\bm{y})-\phi(\bm{x})\geq \langle \bm{\xi}, \,\bm{y}-\bm{x}\rangle, ~\forall \bm{y}\in\mathbb{R}^{d}\}$ the subdifferential of $\phi$ at $\bm{x}$.
The optimal value and optimal solution set of problem~\eqref{1.1} are denoted by $f^{\star}$ and $X^{\star}$, respectively.
We  abbreviate ``almost surely'' as ``a.s''. 
For a collection $\mathcal{G}$ of random variables, $\mathbb{E}[\cdot | \mathcal{G}]$ denotes the conditional expectation. 
Given any $D\subseteq \mathbb{R}^d$,
we denote the distance from $\bm{x}\in\mathbb{R}^d$ to $D$ by $\dist(\bm{x},D) =\inf\{ \|\bm{y}-\bm{x}\| \mid \bm{y}\in D\} $. 
If $D$ is closed and convex, the projection is denoted by $\Pi_{D}(\bm{x})=\arg\min \{ \|\bm{y}-\bm{x}\|\mid \bm{y}\in D\}$. 
We denote $(t)_+ = \max\{t, 0\}$ for $t\in\mathbb{R}$.

\section{Variance Reduced Random Relaxed Projection Method}\label{sec:VR3PM}
We propose a new RPM for solving problem~\eqref{1.1}, namely the variance reduced random relaxed projection method (\alg) in Algorithm~\ref{algoritnm1}. From now on, given any convex function $\phi$, any vector $\bm{x}\in\mathbb{R}^d$ and any subgradient $ \bm{\xi}\in\partial \phi(\bm{x})$, we define 
\begin{equation*}
	\label{eq:H}
	\begin{split}
		H(\phi;\, \bm{x};\, \bm{\xi}) 
		= 
		\begin{cases}
			\{\bm{y}\in \mathbb{R}^d\mid \phi(\bm{x}) + \langle \bm{\xi},\,\bm{y}-\bm{x} \rangle\leq 0\}  & \text{if } \bm{\xi} \neq \bm{0},\\
			\mathbb{R}^d &\text{if } \bm{\xi} = \bm{0}.
		\end{cases}
	\end{split}
\end{equation*}
Note that $H(\phi;\, \bm{x};\, \bm{\xi})$ is a half-space if $\bm{\xi} \neq \bm{0}$. 
\begin{algorithm}[!t]
	\caption{Variance Reduced Random Relaxed Projection Method (\alg)} \label{algoritnm1}
	\begin{algorithmic}[1] 
		\REQUIRE Initial point $\bm{x}^{0}\in C_0$, integers $b\geq 1$ and $r\geq 2$, and a positive sequence $\{\alpha_k\}$. 
		\FOR { $l=0,1,2,\cdots$ }
		\STATE {Set $\tilde{\bm{x}}^{l}=\bm{x}^{lr}$.}
		\FOR{ $s=1,\cdots,r$ }
		\STATE{Set $k=lr+s-1$. Generate i.i.d. uniform indices $I_k=\{i_{k1},\cdots,i_{kb}\}\subseteq[n]$ and compute
			\begin{equation*}
				\bm{v}^{k}=\frac{1}{b}\sum_{i\in I_k}
				(\nabla f_{i}(\bm{x}^{k})-\nabla f_{i}(\tilde{\bm{x}}^{l})) + \nabla f(\tilde{\bm{x}}^{l}).
		\end{equation*} }
		
		\STATE{Generate a random index $j_k\in [m]$, compute a subgradient $\bm{\xi}^{k}\in  \partial \phi_{j_k}(\bm{x}^k)$ and update the iterate
			\[
			\bm{x}^{k+1} =\Pi_{C_0}(\Pi_{H_k}\big(\bm{x}^{k}- \alpha_k\bm{v}^{k}\big)),\]
			where $H_k = H(\phi_{j_k};\, \bm{x}^k;\,\bm{\xi}^k)$.
		}
		\ENDFOR
		\ENDFOR    
	\end{algorithmic} 
\end{algorithm}

Some remarks about \alg are in order. 
First, by definition, $C_{j_k}\subseteq H_k$. Therefore, $H_k$ is an outer approximation of $C_{j_k}$. This explains why we call our method a random ``relaxed'' projection method. 
Second, because of the relaxation in the projection step, the per iteration cost of our method is lower than that of \cite{wang2016stochastic}. Indeed, if $\bm{\xi}^k \neq \bm{0}$, the projection can be computed in closed-form by Lemma~\ref{halfspace}:
\begin{equation*}
	\label{eq:proj_to_H}
	\Pi_{H_k}\big(\bm{x}^{k}- \alpha_k\bm{v}^{k}\big) = \bm{x}^{k}-\alpha_k \bm{v}^{k}-  \frac{ \left( \phi_{j_k}(\bm{x}^k)-\alpha_k\langle \bm{\xi}^k, \bm{v}^{k}\rangle \right)_+ }{\|\bm{\xi}^k\|^2} \bm{\xi}^k.
\end{equation*}
Third, the vector $\bm{v}^k$ is the so-called SVRG gradient estimator~\cite{AVS21}. It is used to reduce the variance of the algorithm and thus improves the convergence speed.
Fourth, the choice of the distribution of the random index $j_k$ is flexible. As long as all the constraints have a positive probability to be chosen~(see Assumption~\ref{A3}), our theoretical framework applies.
Finally, to reduce the variance of \alg due to the random sampling of constraints, we could reformulate the feasible region of problem~\eqref{1.1} by grouping its constraints. Specifically, let $\bar{b} > 0$ be an integer that divides $m$ and define $\bar{m} = m/\bar{b}$. Then, we could re-write the feasible region of problem~\eqref{1.1} as $ C = C_0 \cap ( \bigcap_{ t \in [\bar{m}]} \bar{C}_t )  $,
where 
$\bar{C}_t = \{\bm{x}\in \mathbb{R}^d \mid \max_{j = (t-1)\bar{b} +1,\, \dots,\, t\bar{b} } \phi_j (\bm{x}) \le 0\}$.
This technique can be seen as the constraint analogue of the mini-batch gradient estimator. In our experiments, despite the slightly more complicated subgradient, the overall speed and accuracy could be improved by a suitable grouping.
We should point out that the RPM in \cite{Nedic} is amenable to this technique as it is subgradient-based. However, it might not be worth applying the technique to the RPM in \cite{wang2016stochastic}, as doing so requires computing the projection onto the grouped subset $\bar{C}_t$.

\section{Assumptions}\label{sec:assumptions}
To analyze \alg, the following blanket assumptions on problem~\eqref{1.1} are imposed.
\begin{assumption}\label{ass: problem}
	The following hold.
	\begin{enumerate}[label=(\roman*)]
		\item\label{ass: problem-ii} The set $C_0\subseteq \mathbb{R}^d$ is non-empty, closed and convex. Moreover, for $j\in [m]$, the function $\phi_j$ is proper, closed, and convex.
		
		\item\label{ass: problem-iii} The optimal solution set $X^{\star}$ is non-empty.
		
		\item\label{ass: problem-i} 
		For $i\in [n]$, $f_i$ is differentiable and convex, and there exist $L_i, R_i>0$ such that for any $\bm{x},\bm{y}\in  C_0$,
	 \[
\|\nabla f_i(\bm{x})-\nabla f_i(\bm{y})\|\leq L_i\|\bm{x}-\bm{y}\|+ R_i.
	 \]
	\end{enumerate}
\end{assumption}
\noindent Assumptions~\ref{ass: problem}\ref{ass: problem-ii}-\ref{ass: problem-iii} are standard in the literature of  constrained convex optimization.
Assumption~\ref{ass: problem}\ref{ass: problem-i} is weaker than the requirement on $f_i$ in \cite{Nedic}, which assumes that each $f_i$ has a Lipschitz continuous gradient. Indeed, it can be checked that Assumption~\ref{ass: problem}\ref{ass: problem-i} holds if for any $i\in [n]$, either the function $f_i$ or its gradient~$\nabla f_i$ is Lipschitz continuous.
Also, \cite{wang2016stochastic} requires a similar inequality to hold on~$\mathbb{R}^d$, whereas Assumption~\ref{ass: problem}\ref{ass: problem-i} is assumed to hold on the subset $C_0$.
An immediate consequence of Assumption~\ref{ass: problem}\ref{ass: problem-i} is that
for any $\bm{x}, \bm{y}\in C_0$, 
\begin{align}\label{ineq: Lip-cons-f}
	\|\nabla f(\bm{x}) - \nabla f(\bm{y})\| \leq L\|\bm{x} - \bm{y}\|+R,
\end{align}
where $L = \max_{i \in [n]} L_i$ and $R = \max_{i \in [n]} R_i$.

The assumption below ensures that at each iteration, any constraint will be picked with a positive probability.

\begin{assumption}\label{A3}
	There exists a constant $\rho \in (0,1]$ such that for any $j\in[m]$, 
	\begin{align*}
		\inf_{k\geq 0} P(j_k=j|\mathcal{F}_k)\geq \frac{\rho}{m} \,\,~a.s.,
	\end{align*}
	where $	\mathcal{F}_{k} = \{j_0,\cdots,j_{k-1},i_{01},\cdots,i_{(k-1)b},\bm{x}^0\}$ for $k\ge 1$ and $\mathcal{F}_0 = \{\bm{x}^0\}$.
\end{assumption}

The next assumption concerns the geometry of the feasible region.
\begin{assumption} \label{ass: Line-Regu}
	There exists $\kappa >0$ such that for any $\bm{x}\in C_0$,
	\begin{equation} \label{LR-condition}
		\begin{split}
			\,\dist(\bm{x},C) 
			\leq \, \kappa \max_{{j}\in [m]}\min_{\bm{\xi}_j \in \partial \phi_{j} (\bm{x})}\dist\left(\bm{x},\,  H(\phi_{j};\, \bm{x};\,\bm{\xi}_j) \right).
		\end{split}
	\end{equation}
\end{assumption}

The proposition below asserts that Assumption~\ref{ass: Line-Regu} holds under mild conditions. The proof can be found in Appendix~\ref{appendix-A}. We should emphasize that case~\ref{LR-sufficient-poly} is not new, see~\cite{hoffman2003approximate} for example. The novelty of the proposition lies in case~\ref{LR-sufficient-bounded}.

\begin{proposition} \label{prop:LR-holds}
	Suppose that Assumption~\ref{ass: problem} holds.
	Then, Assumption~\ref{ass: Line-Regu} holds if 
	\begin{enumerate}[label=(\roman*)]
		\item\label{LR-sufficient-poly} $C_0=\mathbb{R}^d$ and $\phi_j$ is affine for all $j\in [m]$, or
		\item\label{LR-sufficient-bounded}  
		the subgradients of $\phi_1,\dots,\phi_m$ have a uniform bound over $C_0$ and $C_0\cap \{ \bm{x}\in \mathbb{R}^d \mid \phi_j(\bm{x}) < 0,\ j \in [m] \} $ is non-empty.
	\end{enumerate}
\end{proposition}

\begin{example} \label{exam:err_bound}
Both optimization problems~\eqref{opt:constr_LASSO} and~\eqref{opt:DRLR} in Examples~\ref{exam:constr_LASSO} and~\ref{exam:DRC}, respectively, satisfy Assumption~\ref{ass: Line-Regu}.
\begin{enumerate}[label=(\roman*)]
\item Since problem~\eqref{opt:constr_LASSO} has only linear constraints, by Proposition~\ref{prop:LR-holds}\ref{LR-sufficient-poly}, it satisfies Assumption~\ref{ass: Line-Regu}.

\item\label{exam:err_bound_DRC} Recall that in problem~\eqref{opt:DRLR}, $d= \ell+1 +n$ and $m=n$. Also, $C_0 = C_{\text{SO}}\times \mathbb{R}_+^n $, where $C_{\text{SO}} = \{ (\bm{u}, \lambda)\in \mathbb{R}^{\ell+1}: \|\bm{u}\|_2\le\lambda \}$ is the second-order cone, and for all $j\in[n]$, $\phi_j (\bm{x}) = \phi(\bm{u},\lambda, \bm{s}) = y_j\,\bm{u}^{\top}\bm{w}_j - \lambda  - s_j$ is linear. Using \cite[Proposition~3 and Corollary~3(c)]{Pang}, there exists a constant $\kappa_1>0$ such that for any $\bm{x} = (\bm{u}^\top, \lambda, \bm{s}^\top)\in C_0$,
\begin{align*}
&\, \dist(\bm{x},C)  = \dist(\bm{x}, (C_{\text{SO}}\times \mathbb{R}_+^n)\cap C_{[m]})  \\
= &\,\dist(\bm{x}, (C_{\text{SO}}\times \mathbb{R}^n)\cap ( \mathbb{R}^{\ell+1} \times \mathbb{R}_+^n)\cap C_{[m]}) \\
\le &\, \kappa_1 \max \lbrace \dist(\bm{x}, ( \mathbb{R}^{\ell+1} \times \mathbb{R}_+^n)\cap C_{[m]} ), ( \|\bm{u}\|-\lambda )_+ \rbrace\\
= &\, \kappa_1 \dist(\bm{x}, ( \mathbb{R}^{\ell+1} \times \mathbb{R}_+^n)\cap C_{[m]} ),
\end{align*}
where the last line follows from $\bm{x}\in C_0 $. Since each $\phi_j$ is linear, the Hoffman error bound~\cite{hoffman2003approximate} and Lemma~\ref{halfspace} imply the existence of a constant $\kappa_2 >0$ such that 
\begin{align*}
&\, \dist(\bm{x}, ( \mathbb{R}^{\ell+1} \times \mathbb{R}_+^n)\cap C_{[m]} )\\
\le&\, \kappa_2 \dist\left(\bm{x},\,  H(\phi_{j};\, \bm{x};\,\bm{\xi}_j) \right).
\end{align*}
Combining the last two display shows that problem~\eqref{opt:DRLR} satisfies Assumption~\ref{ass: Line-Regu}.
\end{enumerate}
\end{example}

Variants of Assumption~\ref{ass: Line-Regu} are utilized to study RPMs.
First, \cite{Nedic} assumes the existence of $\kappa > 0$ such that for any $\bm{x}\in C_0$,
\begin{equation}\label{nedic-err}
	\dist (\bm{x}, C) \leq \kappa \,\mathbb{E}  [(\phi_{j_k}(\bm{x}))_{+} |\, j_0,\cdots,j_{k-1},\bm{x}^0].
\end{equation}
If there exists an upper bound $M>0$ on the sub-differentials $\partial \phi_j (\bm{x})$ that is uniform in $j$ and $\bm{x}$, then condition~\eqref{nedic-err} implies Assumption~\ref{ass: Line-Regu}. To see this, we consider the two cases $\bm{0}\not\in \partial \phi_j(\bm{x})$ and $\bm{0}\in \partial \phi_j(\bm{x})$ separately. If $\bm{0}\not\in \partial \phi_j(\bm{x})$, by Lemma~\ref{halfspace},
\begin{equation*}
	\begin{split}
		&\,(\phi_{j_k}(\bm{x}))_{+}  \le  M  \min_{ \bm{\xi}_{j_k} \in \partial \phi_{j_k}(\bm{x}) } \frac{(\phi_{j_k}(\bm{x}))_{+} }{\|\bm{\xi}_{j_k}\|} \\
		\le &\,  M\max_{j\in[m]} \min_{ \bm{\xi}_j \in \partial \phi_j(\bm{x}) } \dist\left(\bm{x},\,  H(\phi_{j};\, \bm{x};\,\bm{\xi}_j) \right).
	\end{split}
\end{equation*}
If $\bm{0} \in \partial \phi_{j_k}(\bm{x})$, $\bm{x}$ is a minimizer of $\phi_{j_k}$. The non-emptiness of $C_{j_k}$ implies that $\phi_{j_k}(\bm{x}) \le 0$ and hence that $(\phi_{j_k}(\bm{x}))_{+} = 0$. 
Thus, in both cases, condition~\eqref{nedic-err} implies Assumption~\ref{ass: Line-Regu}.

In fact, Assumption~\ref{ass: Line-Regu} is strictly weaker than~\eqref{nedic-err}. 
Consider $\phi_j (\bm{x}) = \bm{x}^{\top} B_j \bm{x}$ for positive definite matrices $B_j$. Then, $ C = \{\bm{0}\} $. If $\bm{x} = \bm{0}$, Assumption~\ref{ass: Line-Regu} holds trivially. If $\bm{x} \neq \bm{0}$, by Lemma~\ref{halfspace}, for any $j\in[m]$, 
\begin{align*}
	&\min_{ \bm{\xi}_j \in \partial \phi_j(\bm{x}) } \dist\left(\bm{x},\,  H(\phi_{j};\, \bm{x};\,\bm{\xi}_j) \right) = \frac{(\phi_j(\bm{x}))_{+}}{\|\nabla \phi_{j}(\bm{x})\|} \\=& \frac{\bm{x}^{\top} B_j \bm{x}}{2\|B_j \bm{x}\|}
	\ge  \frac{\lambda_{\min}(B_j) \|\bm{x}\|^2 }{2\lambda_{\max}(B_j)\|\bm{x}\|} 
	=  \frac{\lambda_{\min}(B_j) }{2\lambda_{\max}(B_j)}  \dist(\bm{x}, C),
\end{align*}
where $\lambda_{\max}(\cdot)$ and $\lambda_{\min}(\cdot)$ denote the maximum and minimum eigenvalues, respectively.
This inequality shows that Assumption~\ref{ass: Line-Regu} holds with $\kappa =  2 \max_{j\in [m]}  \lambda_{\max}(B_j) /\lambda_{\min}(B_j)$.
Also, 
\begin{align*}
	& \frac{\mathbb{E}  [(\phi_{j_k}(\bm{x}))_{+} |\,j_0,\cdots,j_{k-1},\bm{x}^0]}{ \dist(\bm{x}, C) } \le \frac{ \max_{j\in [m]} (\phi_{j}(\bm{x}))_{+} }{ \|\bm{x}\| } \\ 
	=& \max_{j\in[m]} \frac{ \bm{x}^{\top} B_j \bm{x} }{\|\bm{x}\|} \le \|\bm{x}\| \max_{j\in [m]} \lambda_{\max} (B_j),
\end{align*}
where the upper bound vanishes as $\bm{x} \to \bm{0}$, showing that condition~\eqref{nedic-err} fails to hold.

Assumption~\ref{ass: Line-Regu} should also be compared to the classical notion of bounded linear regularity~\cite{bauschke1996projection}:
\[ \dist (\bm{x}, C)  \leq \kappa \max_{j\in [m]} \dist(\bm{x}, C_j) , \]
which has been used for analyzing RPMs in~\cite{wang2016stochastic, Gubin1967TheMO}. Since $C_j \subseteq H(\phi_{j};\, \bm{x};\,\bm{\xi}_j)$, one readily sees that Assumption~\ref{ass: Line-Regu} is stronger than the bounded linear regularity. That a stronger assumption is required by our algorithm is expected as it relies on a weaker projection.

Finally, we remark that Assumption~\ref{ass: Line-Regu} can also be seen as a generalization of the seminal Hoffman error bound~\cite{hoffman2003approximate}, which asserts that the distance of any point to a linear system is linearly bounded by its violation of the linear constraints.
In fact, under case~\ref{LR-sufficient-poly} of Proposition~\ref{prop:LR-holds}, inequality~\eqref{LR-condition} reduces to the Hoffman error bound. Error bound conditions are important subjects in optimization and frequently utilized to study optimization algorithms. For example, as discussed above, conditions similar to (but different from)~\eqref{LR-condition} are employed to study RPMs in~\cite{Nedic} and~\cite{wang2016stochastic}. Going beyond RPMs, error bound conditions also appeared in the study of first-order methods~\cite{luo1993error, liu2017discrete}, second-order methods~\cite{deng2020efficient, yue2019family, yue2019quadratic}, and even manifold optimization algorithms~\cite{liu2017estimation, liu2020unified}.

\section{Convergence Analysis}\label{subsec2}
We now provide a detailed analysis of the convergence behaviour of \alg. The proofs can be found in Appendices~\ref{subsec:converge}-\ref{subsec:converge_rate_QG}.
The first one shows that \alg converges to an optimal solution of problem~\eqref{1.1} almost surely under Assumptions~\ref{ass: problem}--\ref{ass: Line-Regu} and a suitable choice for the step-size $\alpha_k$.

\begin{theorem}\label{convergence1}
	Suppose that Assumptions~\ref{ass: problem}--\ref{ass: Line-Regu} hold. Let $\{\mu_l\}$ be a positive sequence satisfying  
	$\sum_{l=0}^{\infty}\mu_l=\infty$ and  $\sum_{l=0}^{\infty}\mu_l^2< \infty$. 
	Consider Algorithm~\ref{algoritnm1} with $\alpha_k = \mu_l$ for $k=lr +s -1$ with  $l\geq0$ and $s\in [r] $.
	Then, the iterates $\{\bm{x}^{k}\}$ converges almost surely to a point in $X^{\star}$.
\end{theorem}

Our second main theoretical result characterizes the convergence rates of the optimality and feasibility gaps of \alg.
\begin{theorem}\label{rate_pro1}
	Suppose that Assumptions~\ref{ass: problem}--\ref{ass: Line-Regu} hold and that $C_0$ is compact.
	Consider Algorithm~\ref{algoritnm1} with 
	$\alpha_k=\tfrac{\tilde{\alpha}}{\sqrt{k+1}}$, where
	$\tilde{\alpha}\in (0,\tfrac{\rho}{16L m\kappa^2 }]$. 
	Then, for any $K \ge 1$, we have that
	\begin{align*}
		&\mathbb{E}\left[\dist^2( \bar{\bm{x}}^{K},C) \right]\leq \mathcal{O}\left( \tfrac{\log (K)}{K} \right)\\  \text{ and }& \mathbb{E}\left[f(\bar{\bm{x}}^{K})-f^{\star}\right]\leq \mathcal{O} \left(\tfrac{1}{\sqrt{K}} \right),
	\end{align*}
	where $\bar{\bm{x}}^{K}  = \frac{1}{K}\sum_{k=0}^{K-1} \bm{x}^{k}$.
\end{theorem}

With a constant, more conservative step-size, we can improve the bound on the feasibility gap by a factor of $\log K$.
\begin{theorem}\label{rate_convex_constant_stepsize}
	Suppose that Assumptions~\ref{ass: problem}--\ref{ass: Line-Regu} hold and that $C_0$ is compact.
	Consider Algorithm~\ref{algoritnm1} with 
	$\alpha_k \equiv \alpha =\tfrac{\tilde{\alpha}}{\sqrt{K+1}}$, where
	$\tilde{\alpha}\in (0,\tfrac{\rho}{16L m\kappa^2 }]$.
Then, for any $K \ge 1$, we have that
	\begin{align*}
		&\mathbb{E}\left[\dist^2( \bar{\bm{x}}^{K},C) \right]\leq \mathcal{O}\left( \tfrac{1}{K} \right)\\
		\text{ and } & \mathbb{E}\left[f(\bar{\bm{x}}^{K})-f^{\star}\right]\leq \mathcal{O} \left(\tfrac{1}{\sqrt{K}} \right),
	\end{align*}
	where $\bar{\bm{x}}^{K}  = \frac{1}{K}\sum_{k=0}^{K-1} \bm{x}^{k}$.
\end{theorem}

The rate $\mathcal{O} (\tfrac{1}{\sqrt{K}} )$ for the optimality gap in Theorems~\ref{rate_pro1} and~\ref{rate_convex_constant_stepsize} matches the best rate of other RPMs in the literature under the similar assumption setting~\cite{wang2015random, wang2016stochastic, nedic2019, Necoara2022}. As mentioned in the introduction, several works~\cite{ roux2012stochastic, AVR15, AVS21} have successfully accelerated the convergence rate of SGMs via various variance reduction techniques. These works all require that $\nabla f$ is Lipschitz continuous. By making the same assumption, we can also obtain an improved rate for our RPM.

\begin{theorem}
\label{thm:var_redu_rate}
Suppose that Assumptions~\ref{ass: problem}--\ref{ass: Line-Regu} hold with $R_i = 0$ for all $i\in [n]$, that $C$ is compact and that $C_0 = \mathbb{R}^d$. Consider Algorithm~\ref{algoritnm1} with $\alpha_k = \tfrac{\tilde{\alpha}}{(k+1)^{1/3}}$, where $\tilde{\alpha}\in (0,\tfrac{\rho}{40 \sqrt{r}L m\kappa^2 }]$. Then, for any $K\ge 1$, we have that
\begin{align*}
		&\mathbb{E}\left[\dist^2( \bar{\bm{x}}^{K},C) \right]\leq \mathcal{O}\left( \tfrac{1}{K^{2/3}} \right)\\
		\text{ and } & \mathbb{E}\left[f(\bar{\bm{x}}^{K})-f^{\star}\right]\leq \mathcal{O} \left(\tfrac{1}{K^{2/3}} \right),
	\end{align*}
	where $\bar{\bm{x}}^{K}  = \frac{1}{K}\sum_{k=0}^{K-1} \bm{x}^{k}$.
\end{theorem}
\noindent 
The rate $\mathcal{O}( \tfrac{1}{K^{2/3}} )$ for the optimality gap is faster than previous RPMs under comparable assumptions and is achieved by a new analysis framework.
We should point out that the rate $\mathcal{O}( \tfrac{1}{K^{2/3}} )$ for the feasibility gap is slower than those in Theorems~\ref{rate_pro1} and~\ref{rate_convex_constant_stepsize}.

The last result relies on the following assumption, which is called the quadratic growth condition and strictly weaker than strong convexity.

\begin{assumption}
	\label{ass:QG}
	There exists a constant $\nu> 0$ such that
	$f(\bm{x}) - f^{\star} \geq  \frac{\nu }{2}  \dist^2(\bm{x}, X^{\star})$ for any $\bm{x}\in C$.
\end{assumption}

\begin{example}
The feasible region of the constrained least squares problem~\eqref{opt:constr_LASSO} is a polytope $C= \{\bm{x}: B \bm{x} \le \bm{b}\} $ and its objective function $f(\bm{x}) = \| A \bm{x} - \bm{a} \|^2 $ is convex and quadratic over the polytope $C$. By \cite[Proposition~8]{bolte2017error} or \cite[Theorem~8]{necoara2019linear}, $f$ satisfies the quadratic growth condition over $C$. However, if $\ell<d$, then $f$ is not strongly convex.
\end{example}

Stronger rates can be obtained under the quadratic growth condition.
\begin{theorem}\label{rate-qg}
	Suppose that Assumptions \ref{ass: problem}--\ref{ass:QG} hold. Consider Algorithm~\ref{algoritnm1} with     $\alpha_k =  \frac{8}{\nu (l+1)} $ for $k=lr +s-1$ with $l\geq 0$ and $s\in [r]$. Then, for any $K \ge 1$, we have that 
	\begin{align*}
		& \mathbb{E}\left[\dist^2( \bar{\bm{x}}^{K},C)\right]\leq \mathcal{O}\left( \tfrac{1}{K^2} \right) \\
		\text{ and } & \mathbb{E}[f(\bar{\bm{x}}^{K})-f^{\star}]\leq \mathcal{O} \left(\tfrac{1}{K} \right), 
	\end{align*}
	where
	$
	\bar{\bm{x}}^{K} = \frac{3}{K(K+1)(K+2)}\sum_{k\in [K]} k(k+1)\bm{x}^{k}$.
\end{theorem}

\noindent A similar result has been proved for the RPM in~\cite{Necoara2022}.
However, the bounds in~\cite{Necoara2022} hold only if $K$ is sufficiently large, whereas ours hold for any $K\ge 1$. Moreover, \cite{Necoara2022} requires the quadratic growth condition to hold on the larger set $C_0$ but not only~$C$.

\section{Numerical Experiments}\label{sec:experiment}
We then study the empirical performance of \alg through numerical experiments. All the experiments are performed using MATLAB on a PC with Intel Core i5-1135G7 CPU (2.40 GHz). Because of its high accuracy, the optimal solution and optimal value computed by using YALMIP are taken as the ``true'' optimal solution $\bm{x}^{\star}$.

\subsection{Importance of Variance Reduction}\label{sec:exp:SVRG}
In this experiment, we highlight the importance of variance reduction to our algorithm \alg by empirically showing that the SVRG gradient estimator does substantially improve the practical convergence behaviour upon the vanilla gradient estimators. Specifically, we consider the following quadratically constrained quadratic programming problem (QCQP):
\begin{equation}\label{defqcqp}
	\begin{array}{c@{\quad}l}
		\displaystyle \min_{\bm{x}} &\displaystyle \frac{1}{n} \sum_{i\in [n]}  \bm{x}^\top A_i^\top A_i \bm{x} + \bm{a}_i^{\top}\bm{x} \\
		\noalign{\smallskip}
		\mbox{s.t.} & \bm{x}^\top B_j^\top B_j \bm{x} + \bm{b}_j^\top \bm{x} \leq w_j,\quad   j\in [m], \\
		&  \bm{x}\in C_0.
	\end{array}
\end{equation}
Here,
$C_0=[-10, 10]^d$.
The pairs $\{(A_i, \bm{a}_i)\}_{i\in [n]}$ are generated as follows. For each $i\in [n]$, we first generate a random matrix $\tilde{A}_i \in \mathbb{R}^{(p+1)\times d}$ with i.i.d. standard Gaussian entries. Then, the matrix $A_i \in \mathbb{R}^{p\times d}$ and the vector $\bm{a}_i\in \mathbb{R}^d$ are defined as sub-matrices of the normalized matrix $\tilde{A}_i/\norm{\tilde{A}_i}_2 = (A_i^\top \bm{a}_i)^\top$,
where $\|\cdot\|_2$ denotes the operator norm, \ie, the largest singular value.
The pairs $\{(B_j, b_j)\}_{j\in [m]}$ are generated in the same manner.
The constants $w_1,\dots,w_m$ are i.i.d. uniform random variables on $[0, 0.5]$.
We compare \alg with three variants of random relaxed projection methods (R$^2$PMs) obtained by replacing the SVRG gradient estimator $\bm{v}^{k}$ in Algorithm~\ref{algoritnm1} with the standard gradient estimator using a single summand, the mini-batch gradient estimator using $b$ summands and the full gradient using all $n$ summands. These three R$^2$PMs are denoted, respectively, as R$^2$PM-$1$, R$^2$PM-$b$ and R$^2$PM-$n$. The step-sizes for all algorithms are chosen to scale as $\alpha_k = \mathcal{O}(k^{-0.51})$ with optimally tuned constants. The constraint grouping technique is applied to all tested algorithms. Specifically, we group the constraints into $\bar{m}= m/10$ groups of $\bar{b}=10$ constraints.

\begin{figure}[!t]
	\centering
	\begin{subfigure}[t]{1\columnwidth }
		\centering
		\includegraphics[width=0.48\textwidth]{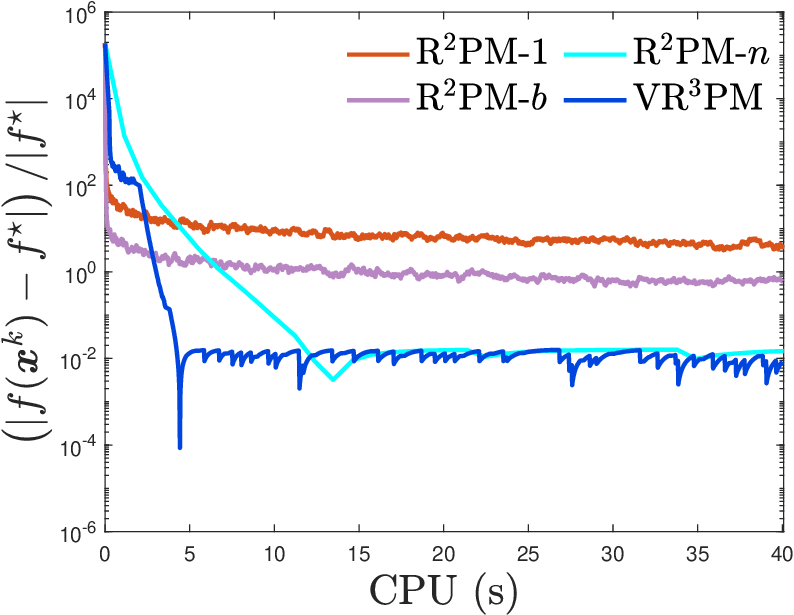}
		\includegraphics[width=0.48\textwidth]{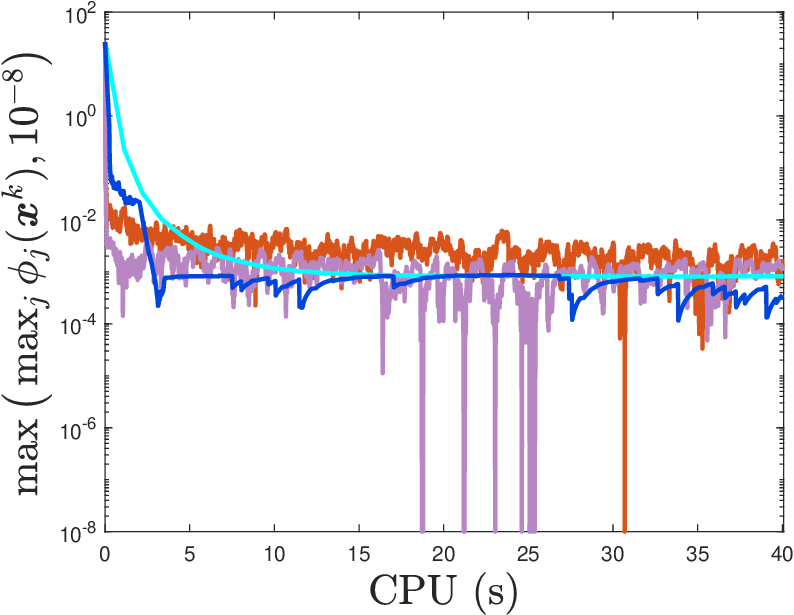}
		\caption{$(n,m,d, p, q)=(3000, 3000, 200, 200,200)$\vspace{4mm}}
		\label{QCQP-VR-compare-a}
	\end{subfigure}
	\begin{subfigure}[t]{1\columnwidth}
		\centering
		\includegraphics[width=0.48\textwidth]{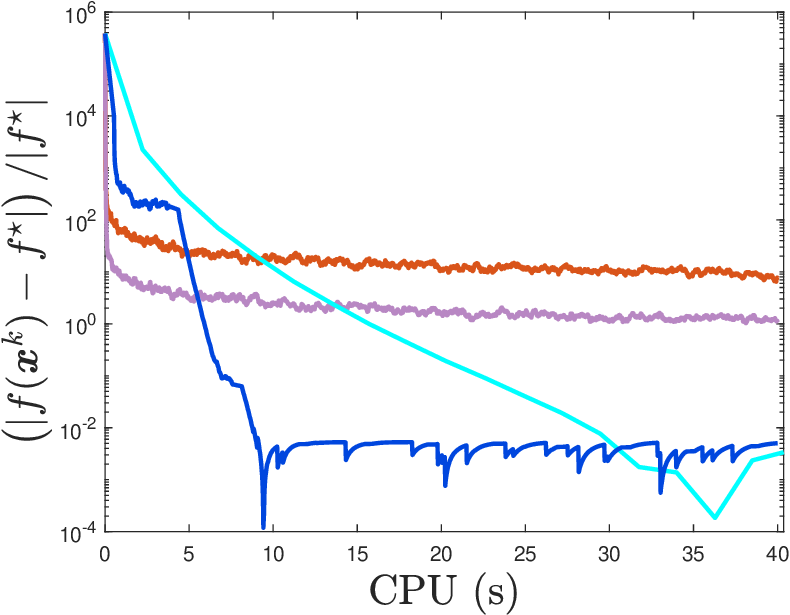}
		\includegraphics[width=0.48\textwidth]{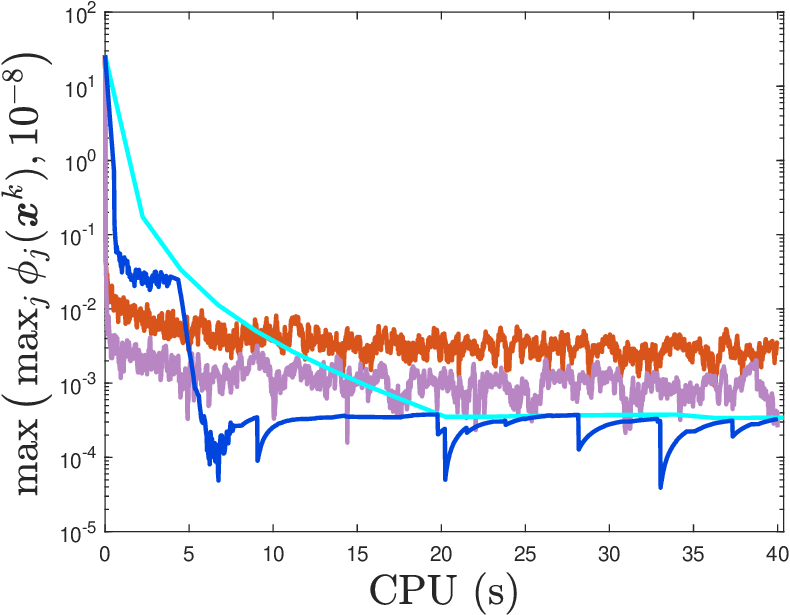}
		\caption{$(n,m,d, p,q)=(6000, 6000, 200, 200, 200)$}
		\label{QCQP-VR-compare-b}
	\end{subfigure}
	\caption{Comparison of VR$^{3}$PM and other RPMs using vanilla gradient estimators on problem~\eqref{defqcqp}.}
	\label{QCQP-VR-compare}
\end{figure}

The results on problem~\eqref{defqcqp} with parameters $(n, m ,d, p, q)=(3000, 3000, 200, 200, 200)$ and $(n, m ,d, p, q)=(6000, 6000, 200, 200, 200)$  are shown in Figure~\ref{QCQP-VR-compare}, where the left and right panels show the optimality gap and constraint violation against the CPU time (in second), respectively.
We can see from Figure~\ref{QCQP-VR-compare} that in terms of objective value, our algorithm \alg is faster and reaches a higher accuracy than the other three R$^2$PMs. As for the constraint violation, our algorithm \alg performs on par with R$^2$PM-$1$ and R$^2$PM-$b$ and outperforms the full gradient variant R$^2$PM-$n$.

\subsection{Comparison with the RPMs in \cite{Necoara2022} and \cite{wang2016stochastic}}\label{sec:exp:comparison}
Our second experiment aims at comparing the performance of \alg with the existing RPMs in~\cite{Necoara2022} by Necoara and Singh and in~\cite{wang2016stochastic} by Wang and Bertsekas, denoted as RPM-NS, and RPM-WB, respectively.
RPM-WB is not amenable to the constraint grouping technique. So, we apply the technique only to our algorithm and RPM-NS. The constraint group size and step-sizes are chosen similarly as in Section~\ref{sec:exp:SVRG}.

\begin{figure}[!t]
	\centering
	\begin{subfigure}[t]{1\columnwidth }
		\centering
		\includegraphics[width=0.48\textwidth]{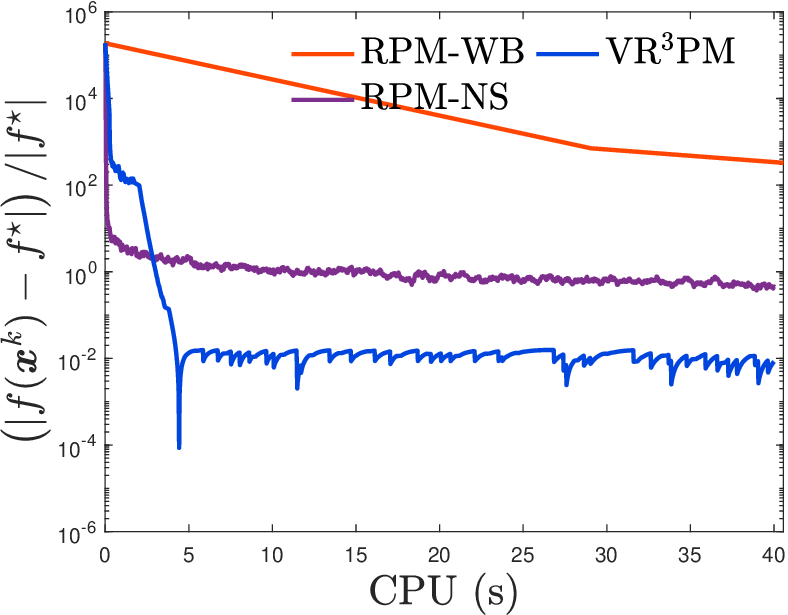}
		\includegraphics[width=0.48\textwidth]{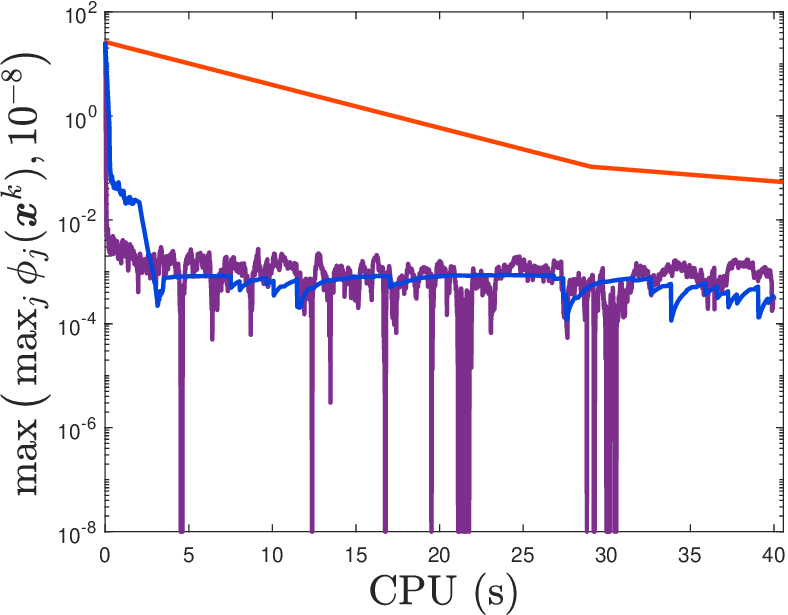}
		\caption{$(n,m,d, p,q)=(3000, 3000, 200, 200,200)$\vspace{4mm}}
		\label{QCQP-all-compare-a}
	\end{subfigure}
	\begin{subfigure}[t]{1\columnwidth }
		\centering
		\includegraphics[width=0.48\textwidth]{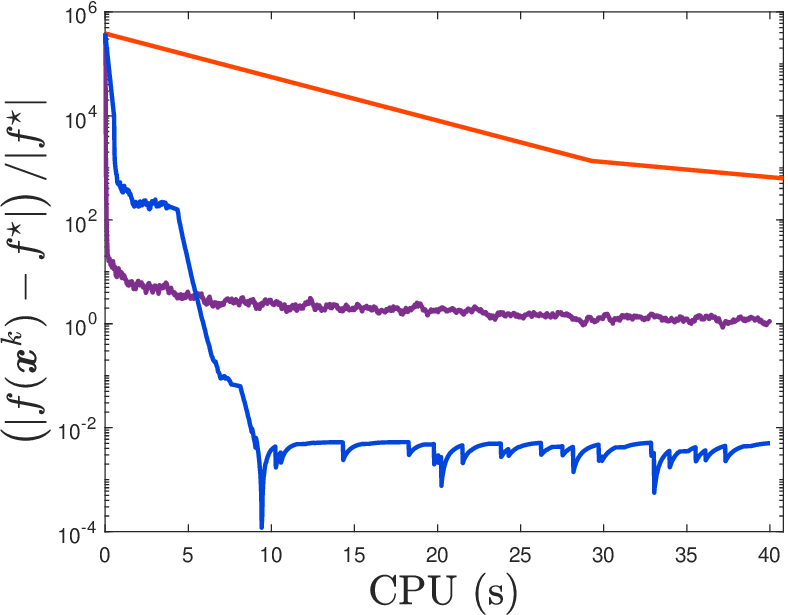}
		\includegraphics[width=0.48\textwidth]{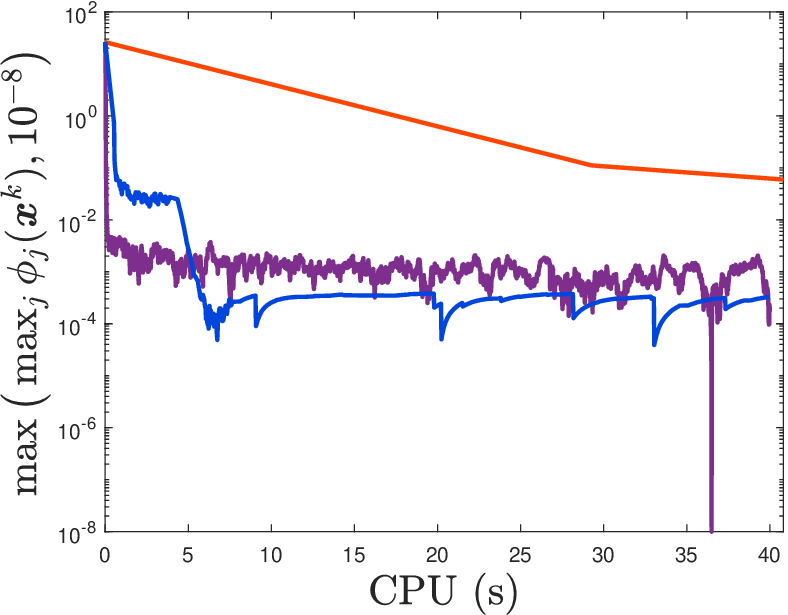}
		\caption{$(n,m,d, p,q)=(6000, 6000, 200,200,200)$}
		\label{QCQP-all-compare-b}
	\end{subfigure}
	\caption{Comparison of VR$^{3}$PM and the RPMs in~\cite{Necoara2022} and~\cite{wang2016stochastic} on problem~\eqref{defqcqp}.}
	\label{QCQP-all-compare}
\end{figure}

The results on problem~\eqref{defqcqp} with the same setting and parameters as in Section~\ref{sec:exp:SVRG} are shown in Figure~\ref{QCQP-all-compare}, where the left and right panels show the optimality gap and constraint violation against the CPU time (in second), respectively.
From Figure~\ref{QCQP-all-compare}, we can see that our algorithm \alg performs substantially better than the competing algorithms RPM-NS and RPM-WB in terms of objective value. As for constraint violation, our algorithm and RPM-NS perform on par, but both better than RPM-WB.

\subsection{Applications to Downlink Beamforming}
Our next experiment concerns the downlink beamforming problem in wireless communication~\cite{wu2017sdr, Luo2006AnIT, Wiesel2006LinearPV}. Specifically, we consider a single base station equipped $d$ antennas, transmitting data stream to $n$ users. 
For each user $i\in [n]$, the signal received is given by
\[ y_i = \bm{h}_i^{\mathsf{H}}\left( \sum_{j\in[n]} \nu_j \bm{x}_j\right) + \omega_i, \]
where $(\cdot)^{\mathsf{H}}$ denotes the Hermitian transpose, and $\bm{h}_i \in \mathbb{C}^d$ models the downlink channel, $\omega_i \in \mathbb{C}$ is the error, $\nu_i$ is the information signal and $\bm{x}_i \in\mathbb{C}^d$ represents the beamformer for user $i$. 
The quality of the signal received at user $i$ can be measured by the signal-to-interference-plus-noise ratio (SINR):
\[ \mathrm{SINR}_i = \frac{ | \bm{h}_i^{\mathsf{H}} \bm{x}_i |^2 }{  \sum_{j\neq i} |\bm{h}_i^{\mathsf{H}} \bm{x}_j|^2 + \sigma^2 }, \]
where $\sigma^2 $ is the variance of the error $\omega_i$. One formulation of the beamforming problem is to minimize the transmission power subject to SINR constraints~\cite{Luo2006AnIT}:
\begin{equation}\label{beamforming-1}
	\begin{array}{c@{\quad}l}
		\displaystyle\min_{\bm{x}_1,\dots, \bm{x}_n} &\displaystyle \sum_{i=1}^n \|\bm{x}_i\|^2 \\
		\noalign{\smallskip}
		\mbox{s.t.} & \mathrm{SINR}_i \geq \gamma_i,\quad  i\in [n],
	\end{array}
\end{equation}
where $\gamma_1,\dots, \gamma_n >0$ are the desired SINRs. Problem~\eqref{beamforming-1} is non-convex since the SINR constraints are non-convex. However, it is shown in~\cite{Wiesel2006LinearPV} that problem~\eqref{beamforming-1} is equivalent to a convex, second-order cone program over real numbers:
\begin{equation}\label{beamforming-SOCP}
	\begin{array}{c@{\quad}l}
		\displaystyle\min_{\hat{\bm{x}}_1,\dots, \hat{\bm{x}}_n} &\displaystyle \sum_{i=1}^n \|\hat{\bm{x}}_i\|^2 \\
		\noalign{\smallskip}
		\mbox{s.t.} & \hat{\bm{h}}_i^\top \hat{\bm{x}}_i \geq \sqrt{\gamma_i \sum_{j\neq i} (\hat{\bm{h}}_j^\top \hat{\bm{x}}_j)^2 + \gamma_i\sigma^2},~  i \in [n],
	\end{array}
\end{equation}
where $\hat{\bm{h}}_i\in \mathbb{R}^{2d}$ and $\hat{\bm{x}}_i \in \mathbb{R}^{2d}$ are real vectors obtained by stacking the real and imaginary parts of $\bm{h}_i$ and $\bm{x}_i$, respectively, for $i\in [n]$.
Using a similar argument as in Example~\ref{exam:err_bound}\ref{exam:err_bound_DRC}, one could also show that problem~\eqref{beamforming-SOCP} satisfies Assumption~\ref{ass: Line-Regu}.

In this experiment, we set $\sigma = 1$ and the SINR threshold for all users to be $-13 $ dB (\ie, $\gamma_i \approx 0.0501$). The vectors $\hat{\bm{h}}_1,\dots, \hat{\bm{h}}_n$ are i.i.d. $2d$-dimensional standard Gaussian random distribution with zero mean and identity covariance.

\begin{figure}[!t]
	\begin{subfigure}[t]{1\columnwidth }
		\centering
		\includegraphics[width=0.48\textwidth]{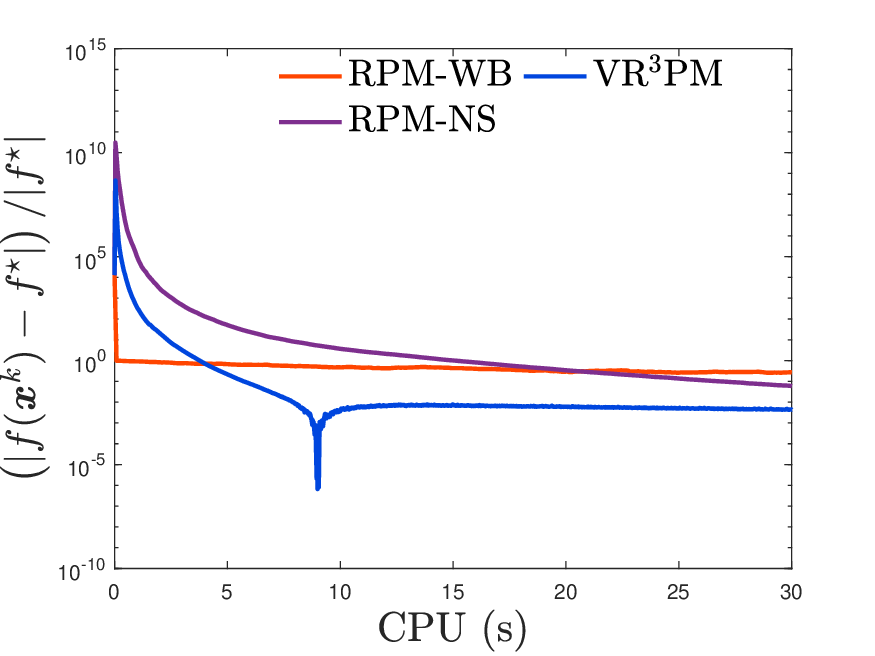}
		\includegraphics[width=0.48\textwidth]{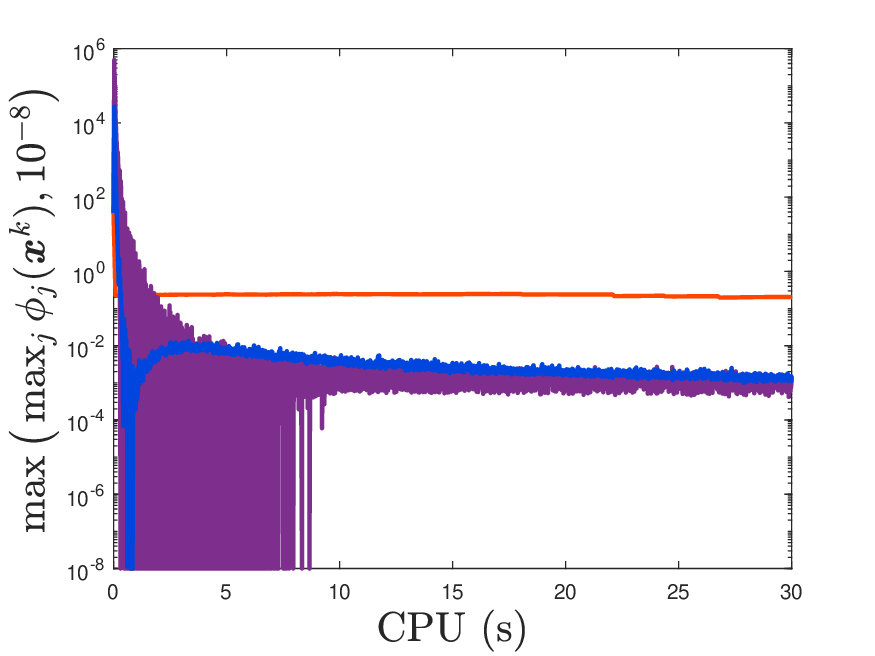}
		\caption{$(n,d)=(200,25)$}
		\label{BF-VR-compare-a}
	\end{subfigure}
	\begin{subfigure}[t]{1\columnwidth }
		\centering
		\includegraphics[width=0.48\textwidth]{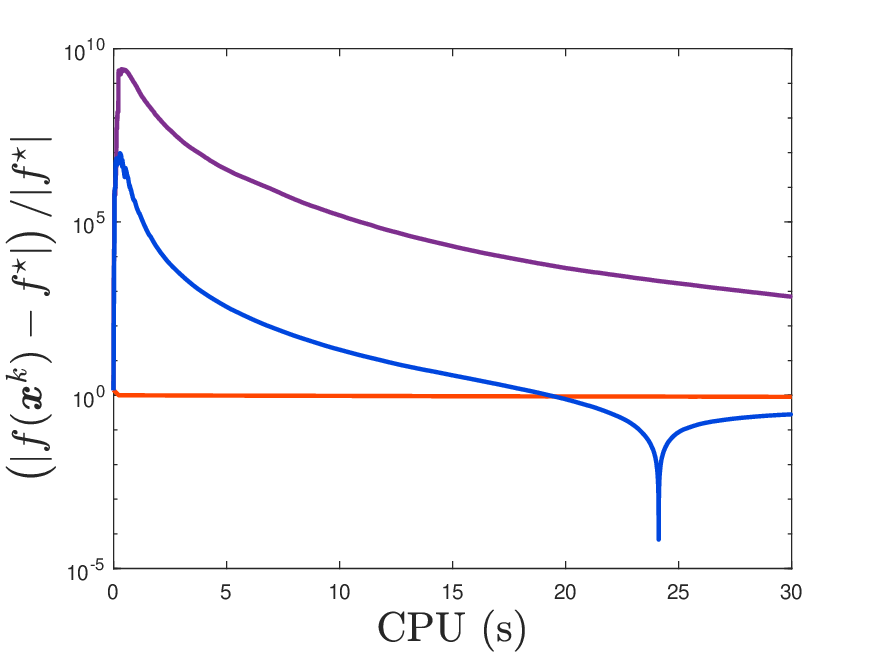}
		\includegraphics[width=0.48\textwidth]{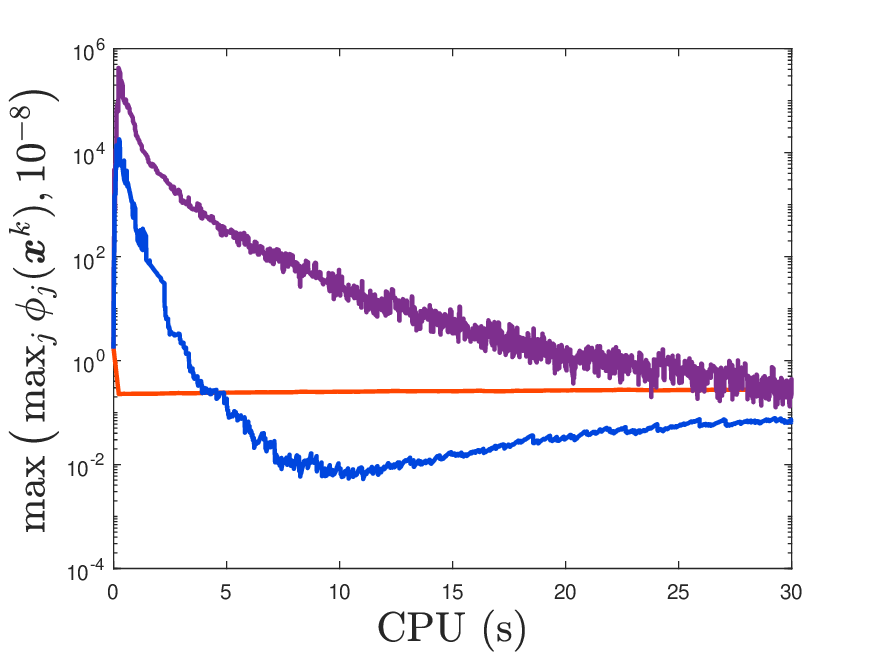}
		\caption{$(n,d)=(600,20)$}
		\label{BF-VR-compare-b}
	\end{subfigure}
	\caption{Comparison of VR$^{3}$PM and the RPMs in~\cite{Necoara2022} and~\cite{wang2016stochastic} on problem~\eqref{beamforming-SOCP}.}
	\label{BF-VR-compare}
\end{figure}

The results on problem~\eqref{beamforming-SOCP} with parameters $(n, d)=(200, 25)$ and $(n,d)=(600, 20)$ are shown in Figure~\ref{BF-VR-compare}, where we can see that \alg performs better than RPM-NS and RPM-WB in terms of both  optimality and feasibility gaps.

\subsection{Applications to Robust Classification}
Finally, we compare our algorithm against the RPM-NS and RPM-WB on the distributionally robust classification problem in Example~\ref{exam:DRC}.
In this experiment, we use real data-sets\footnote{\href{https://www.csie.ntu.edu.tw/~cjlin/libsvmtools/datasets/}{https://www.csie.ntu.edu.tw/~cjlin/libsvmtools/datasets/}} {\bf a6a} with $(n,\ell) = (11220, 122)$ and {\bf a7a} with $(n,\ell) =(16100, 122)$. In both data-sets, the sample size $n$ is very large. This makes problem~\eqref{opt:DRLR} extremely computationally challenging as the number of decision variables, the number of summands in the objective and the number of constraints all increase linearly in $n$. We compare our algorithm \alg against RPM-NS and RPM-WB with all the algorithmic parameters similarly chosen.

\begin{figure}[!t]
	\centering
	\begin{subfigure}[t]{1\columnwidth }
		\centering
		\includegraphics[width=0.48\textwidth]{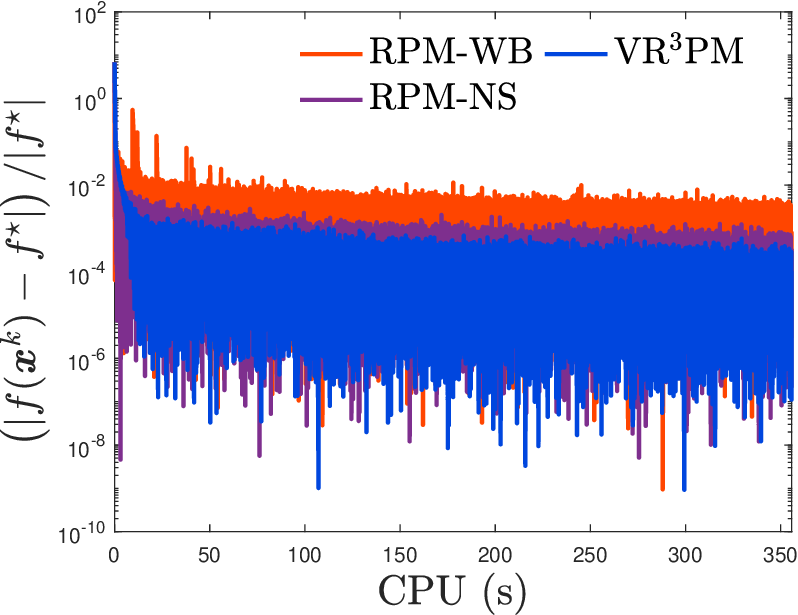}
		\includegraphics[width=0.48\textwidth]{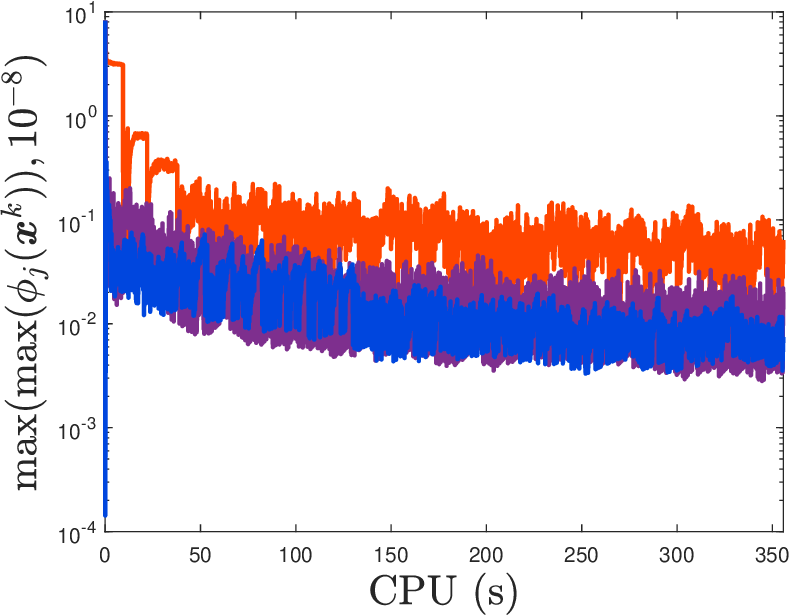}
		\caption{Results on {\bf a6a} with $(n,\ell)=(11220,122)$}
		\label{DRLR-all-compare-b}
	\end{subfigure}
	\begin{subfigure}[t]{1\columnwidth }
		\centering
		\includegraphics[width=0.48\textwidth]{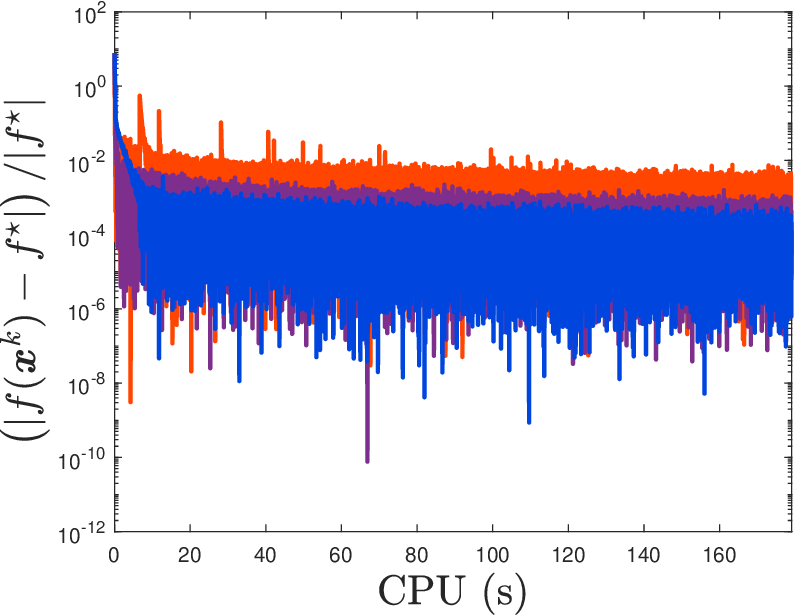}
		\includegraphics[width=0.48\textwidth]{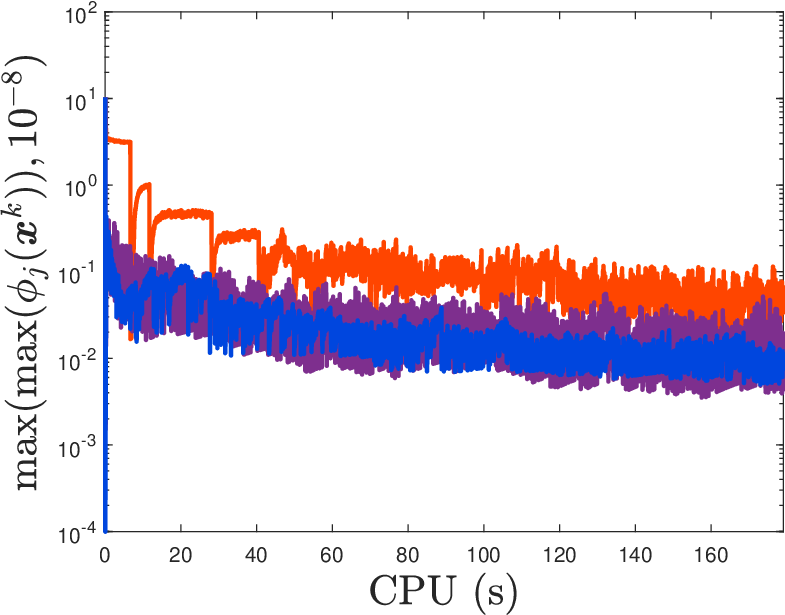}
		\caption{Results on {\bf a7a} with $(n,\ell)=(16100,122)$}
		\label{DRLR-all-compare-c}
	\end{subfigure}
	\caption{Comparison of VR$^{3}$PM and the RPMs in~\cite{Necoara2022} and~\cite{wang2016stochastic} on problem~\eqref{opt:DRLR}.}
	\label{DRLR-all-compare}
\end{figure}

The results are presented in Figure~\ref{DRLR-all-compare}, where we can see that \alg performs better than the two competing RPMs in terms of both the optimality gap and constraint violation.

\appendix

\section{Proof of Main Results} \label{appendix-C}

\subsection{Lemmas}\label{appendix-B}
The following are standard results on projections, see for example~\cite[Lemma~2]{wang2016stochastic} and \cite[Theorem~5.4 and Lemma~6.26]{beck2017first}.

\begin{lemma}\label{lema: Proj-Property}
	Let $D\subseteq \mathbb{R}^d$ be a non-empty, closed and convex set. The following hold.
	\begin{enumerate}[label=(\roman*)]
		\item\label{lema: proj-property-nonexpan} For  any $\bm{x}$, $\bm{y}\in \mathbb{R}^d$,  
		\[
		\|\Pi_{D}(\bm{x})-\Pi_{D}(\bm{y})\|\leq \|\bm{x}-\bm{y}\|;
		\]
		\item\label{lema: proj-property-tri-in}  For any $\bm{x}\in \mathbb{R}^d$ and $\bm{y}\in D$,~
		\[ \|\Pi_D(\bm{x})-\bm{y}\|^2\leq \|\bm{x}-\bm{y}\|^2-\|\bm{x}-\Pi_D(\bm{x})\|^2.\]
	\end{enumerate}
\end{lemma}

\begin{lemma}\label{halfspace}
	Let $c\in\mathbb{R}$, $\bm{\zeta}\in\mathbb{R}^d\setminus \{\bm{0}\}$ and $H=\{\bm{x}\mid \bm{\zeta}^{T}\bm{x}\leq c\}$. Then
	\begin{align}\notag
		\Pi_{H}(\bm{u})=\bm{u}-\frac{(\bm{\zeta}^{T}\bm{u}-c)_{+}}{\|\bm{\zeta}\|^2}\bm{\zeta}.
	\end{align}
\end{lemma}

The following lemma are useful to our development, see~\cite{Nedic,wang2016stochastic} for example.
\begin{lemma} \label{lema: supermartingale}
	Let  $\{a^k\}$, $\{u^k\}$, $\{t^k\}$ and $\{d^k\}$ be sequences of non-negative random variables satisfying
	\begin{align*}
		& \mathbb{E}[a^{k+1}| a^0,\cdots,a^k,u^0,\cdots,u^k,t^0,\cdots,t^k,d^0,\cdots,d^k]\\
		\leq &(1+t^k)a^k-u^k+d^k~~\text{for all}~~k\geq 0~a.s.
	\end{align*}
	Suppose that 
	\[ \sum_{k=0}^{\infty}t^k<\infty\quad \text{and} \quad \sum_{k=0}^{\infty}d^k<\infty \qquad a.s. \]
	Then,
	\[ \sum_{k=0}^{\infty}u^k< \infty\,\, a.s., \quad\text{and} \quad \lim_{k\rightarrow 0}a^k=a \qquad a.s.,\]
	for some non-negative random variable $a$.
\end{lemma}

The following lemma bounds the distance of an iterate to the feasible region in terms of its distance to the corresponding half-space at that iteration.

\begin{lemma}\label{delta}
	Suppose that Assumptions~\ref{A3}--\ref{ass: Line-Regu} hold.
	Then, Algorithm~\ref{algoritnm1} satisfies that for any $k\ge 0$,
	\begin{equation*}
		\begin{split}
			\, \|\bm{x}^k - \Pi_C(\bm{x}^k)\|^2 
			\leq \, m\kappa^2 \rho^{-1}\, \mathbb{E}\left[\|\bm{x}^k-\Pi_{H_k}(\bm{x}^k)\|^2 | \mathcal{F}_k \right] \; a.s.  
		\end{split}
	\end{equation*}
\end{lemma}

\begin{proof}
	We have that for any $j\in [m]$,
	\begin{align*}
		& \, \mathbb{E}\left[\|\bm{x}^k - \Pi_{H_k}(\bm{x}^k)\|^2|\mathcal{F}_k\right] \\
		= &\,\sum_{j' \in [m]} P(j_k = j'|\mathcal{F}_k)\|\bm{x}^k-\Pi_{ H(\phi_{j'};\, \bm{x}^k;\,\bm{\xi}_{j'}) }(\bm{x}^k)\|^2\\
		\geq & \,\frac{\rho}{m} \|\bm{x}^k-\Pi_{ H(\phi_{j};\, \bm{x}^k ; \,\bm{\xi}^k))  }(\bm{x}^k)\|^2 \\
		\ge &\,\frac{\rho}{m} \min_{\bm{\xi}_{j} \in \partial \phi_{j} (\bm{x}^k)}\dist^2(\bm{x}^k, H(\phi_{j};\, \bm{x}^k;\,\bm{\xi}_j) ),
	\end{align*}
	where the inequality follows from Assumption~\ref{A3}. Therefore, using Assumption~\ref{ass: Line-Regu}, we obtain
	\begin{align*}
		&  \mathbb{E} \left[\|\bm{x}^k - \Pi_{H_k}(\bm{x}^k)\|^2|\mathcal{F}_k\right] \\
		\ge    & \, \frac{\rho}{m} \max_{j\in [m]} \min_{\bm{\xi}_{j} \in \partial \phi_{j} (\bm{x}^k)}\dist^2(\bm{x}^k, H(\phi_{j};\, \bm{x}^k;\,\bm{\xi}_j) )\\
		\ge  & \, \frac{\rho}{m \kappa^2} \|\bm{x}^k - \Pi_C(\bm{x}^k)\|^2,
	\end{align*}
	which completes the proof.
\end{proof}

The next lemma should be compared with \cite[Assumption~1(c)]{wang2016stochastic} and, to some extent, illustrates why the SVRG gradient estimator cannot be directly used in the algorithmic framework of \cite{wang2016stochastic}.
\begin{lemma}\label{lema: SVRG-type-Property}
	Suppose that Assumption \ref{ass: problem} holds.
	Consider Algorithm~\ref{algoritnm1}. Then, for any $\bm{x}^{\star}\in X^{\star}$ and $k=lr+s-1$, where $l\geq0$, $s\in[r]$, we have
	\begin{align*}
		&\mathbb{E}[\|\bm{v}^{k}\|^2\mid\mathcal{F}_k]  \\
		\leq &\,  8L^2\|\bm{x}^k-\bm{x}^{\star}\|^2+16L^2\|\tilde{\bm{x}}^{l}-\bm{x}^{\star}\|^2 + 12R^2 +4\|\nabla f(\bm{x}^{\star})\|^2.
	\end{align*}
\end{lemma}

\begin{proof}
	We have
	\begin{equation*} 
		\begin{split}
			&\, \mathbb{E}\left[\|\bm{v}^{k}\|^2\mid\mathcal{F}_{k}\right] \\
			= & \,\mathbb{E}\left[ \bigg\|\frac{1}{b}\sum_{i\in I_k}(\nabla f_i(\bm{x}^k)-\nabla f_i(\tilde{\bm{x}}^{l}))+ \nabla f(\tilde{\bm{x}}^{l}) \bigg\|^2 \,\bigg|\, \mathcal{F}_k \right] \\
			\le &\, 2 \mathbb{E}\left[\bigg\|\frac{1}{b}\sum_{i\in I_k}(\nabla f_i(\bm{x}^k) -\nabla f_i(\tilde{\bm{x}}^{l}))\bigg\|^2 \,\bigg|\,\mathcal{F}_k \right] \\
			&+ 2\mathbb{E}\left[\|\nabla f(\tilde{\bm{x}}^{l}))\|^2\mid\mathcal{F}_k \right]\\
			\le & \, 4L^2 \|\bm{x}^k-\tilde{\bm{x}}^{l}\|^2+ 4R^2 +4\|\nabla f(\tilde{\bm{x}}^{l})-\nabla f(\bm{x}^{\star})\|^2 \\
			&+4\|\nabla f(\bm{x}^{\star})\|^2 \\   
			\le & \,
			8L^2\|\bm{x}^k-\bm{x}^{\star}\|^2 + 16L^2\|\tilde{\bm{x}}^{l}-\bm{x}^{\star}\|^2  + 12R^2 +4\|\nabla f(\bm{x}^{\star})\|^2,
		\end{split}
	\end{equation*}
	where the equality follows from  the definition of $\bm{v}^{k}$, the second inequality follows from Assumption~\ref{ass: problem}\ref{ass: problem-i}, $\bm{x}^{k}\in C_0$ and $\tilde{\bm{x}}^{l}\in C_0$, and the last inequality follows from inequality~\eqref{ineq: Lip-cons-f}. This completes the proof.
\end{proof}

The following technical lemma is also used in the proof of our main results.

\begin{lemma}\label{lema: 3.5}
	Suppose that Assumption \ref{ass: problem} holds. 
	Consider Algorithm~\ref{algoritnm1}.  Then, for any $\bm{x}^{\star}\in X^{\star}$, $k\geq 0$ and $\lambda>0$, we have
	\begin{align*}
		& \, 2 \alpha_k \mathbb{E} \left[ \langle \bm{v}^{k}, \bm{x}^k - \bm{x}^{\star} \rangle \mid \mathcal{F}_k \right]  \\
		\geq & -(2 \alpha_k L + \tfrac{4}{\lambda})\|\bm{x}^k-\Pi_C(\bm{x}^k)\|^2  -  \alpha_k^2 \lambda  \left( 2 R^2 + \|\nabla f(\bm{x}^{\star})\|^2 \right)  \\
		&\,   - \alpha_k^2  L^2 \lambda \|\bm{x}^k-\bm{x}^{\star}\|^2 + 2 \alpha_k \left( f( \Pi_C( \bm{x}^k ) ) - f^{\star} \right).
	\end{align*}
\end{lemma}

\begin{proof}
	First, for any $k\geq0$,
	\begin{equation}\label{ineq: lema3.5-1}
		\begin{split}
			&\, 2\alpha_k \mathbb{E}\left[\langle \bm{v}^{k}, \, \bm{x}^k - \bm{x}^{\star} \rangle \mid \mathcal{F}_k\right]  \\
			= &\, 2\alpha_k \langle \nabla f(\bm{x}^k),\,\bm{x}^k-\bm{x}^{\star}\rangle
			\ge  2\alpha_k  (f(\bm{x}^k)-f(\bm{x}^{\star}) )    \\
			= & \, 2\alpha_k \left(f(\bm{x}^k) - f(\Pi_C (\bm{x}^{k}))+f(\Pi_C (\bm{x}^{k})) -f^{\star}\right) \\
			\ge & \,  2\alpha_k   \langle \nabla f(\Pi_C(\bm{x}^k)),\, \bm{x}^k - \Pi_C(\bm{x}^k)\rangle \\
			&\, +2\alpha_k \left( f(\Pi_C(\bm{x}^k))-f^{\star} \right), 
		\end{split}
	\end{equation}
	where the first equality follows from the definitions of $\bm{v}^{k}$, the random index subset $I_k$ and the collection $\mathcal{F}_{k}$, the first inequality from the convexity of $f$, and the second inequality from the convexity of $f$.
	Also, we have
	\begin{equation}\label{ineq:proof:1}
		\begin{split}
			&2 \alpha_k\langle  \nabla f(\Pi_C(\bm{x}^k)),\, \bm{x}^k - \Pi_{C} (\bm{x}^k) \rangle \\
			\geq & -2 \alpha_k \|\nabla f(\Pi_C(\bm{x}^k)) - \nabla f(\bm{x}^k) \| \| \bm{x}^k-\Pi_C(\bm{x}^k) \|   \\
			&  -2 \alpha_k   \| \nabla f(\bm{x}^k) \| \| \bm{x}^k - \Pi_{C}( \bm{x}^k ) \| \\
			\geq & -2 \alpha_k   ( L\|\bm{x}^k-\Pi_C(\bm{x}^k)\|+R)\|\bm{x}^k-\Pi_C(\bm{x}^k ) \|  \\
			&-2\alpha_k\| \nabla f(\bm{x}^k)\|\|\bm{x}^k-\Pi_{C}(\bm{x}^k)\| \\
			\geq & - 2 \alpha_k L\|\bm{x}^k-\Pi_C(\bm{x}^k)\|^2 -2 \alpha_k R\| \bm{x}^k-\Pi_C(\bm{x}^k) \|  \\
			& -2 \alpha_k\| \nabla f(\bm{x}^k)-\nabla f(\bm{x}^{\star})\|\|\bm{x}^k-\Pi_{C}(\bm{x}^k)\| \\ 
			&-2 \alpha_k\|\nabla f(\bm{x}^{\star})\|\|\bm{x}^k-\Pi_{C}(\bm{x}^k)\|,
		\end{split}
	\end{equation}
	where the first inequality follows from the Cauchy-Schwarz inequality, the second from inequality~\eqref{ineq: Lip-cons-f}, $\bm{x}^{k}\in C_0$ and $\Pi_{C}(\bm{x}^{k})\in C\subseteq C_0$, and the third from the triangle inequality. We then bound the second, third and fourth terms on the last line of \eqref{ineq:proof:1}. We will use multiple times the fact that $2|a_1 a_2 |\leq \lambda a_1^2+\frac{1}{\lambda}a_2^2$ for all $a_1, a_2\in \mathbb{R}$ and $\lambda>0$.
	The second term can be bounded as 
	\begin{equation} \label{ineq:proof:2}
		\begin{split}
			& -2 \alpha_k R\| \bm{x}^k-\Pi_C(\bm{x}^k) \| \\
			\geq &\, -  \alpha_k^2 R^2 \lambda - \frac{1}{\lambda} \| \bm{x}^k - \Pi_C(\bm{x}^k) \|^2.  
		\end{split}
	\end{equation}
	The third term can be bounded as
	\begin{equation}
		\label{ineq:proof:3}
		\begin{split}
			& -2 \alpha_k\| \nabla f(\bm{x}^k) - \nabla f(\bm{x}^{\star})\| \| \bm{x}^k - \Pi_{C}(\bm{x}^k) \| \\
			\ge & -2 \alpha_k (L\|\bm{x}^k-\bm{x}^{\star}\|+ R) \|\bm{x}^k-\Pi_{C}(\bm{x}^k)\| \\ 
			\ge & - \alpha_k^2L^2\lambda\|\bm{x}^k-\bm{x}^{\star}\|^2-\frac{1}{\lambda}\|\bm{x}^k-\Pi_{C}(\bm{x}^k)\|^2 \\
			&-2\alpha_k R \|\bm{x}^k - \Pi_C(\bm{x}^k) \|, \\
			\ge &
			- \alpha_k^2L^2\lambda\|\bm{x}^k-\bm{x}^{\star}\|^2-  \alpha_k^2 R^2 \lambda  \\
			&-\frac{2}{\lambda}\|\bm{x}^k-\Pi_{C}(\bm{x}^k)\|^2,
		\end{split}
	\end{equation}
	where the first inequality follows from inequality~\eqref{ineq: Lip-cons-f} and the last from~\eqref{ineq:proof:2}.
	And the fourth term can be bounded as
	\begin{equation}\label{ineq:proof:4}
		\begin{split}
			&    -2 \alpha_k\| \nabla f(\bm{x}^{\star})\| \|\bm{x}^k-\Pi_{C}(\bm{x}^k) \|\\
			\geq & - \alpha_k^2\lambda\| \nabla f(\bm{x}^{\star}) \|^2-\frac{1}{\lambda} \|\bm{x}^k - \Pi_{C}(\bm{x}^k) \|^2.  
		\end{split}
	\end{equation}
	Substituting inequalities~\eqref{ineq:proof:1}, \eqref{ineq:proof:2}, \eqref{ineq:proof:3} and~\eqref{ineq:proof:4} into~\eqref{ineq: lema3.5-1} yields the desired result.
\end{proof}

The next lemma establishes a recursion for the distance to optimality $\|\bm{x}^k - \bm{x}^\star\|$.

\begin{lemma}\label{lemma: one-iteration}
	Suppose that Assumptions \ref{ass: problem}--\ref{ass: Line-Regu} hold. Consider Algorithm~\ref{algoritnm1}.
	Then, for any $\lambda >0$, $\bm{x}^{\star}\in X^{\star}$ and $k=lr+s-1$, where $l\geq0$, $s\in[r]$, 
	we have
	\begin{equation*}
		\begin{split}
			&\, \mathbb{E} \left[\|\bm{x}^{k+1}-\bm{x}^{\star}\|^2\mid\mathcal{F}_{k} \right]  \\
			\leq\, &\, (1 + \alpha_k^2(24 L^2+L^2\lambda))\|\bm{x}^k-\bm{x}^{\star}\|^2  \\
			&\, + \alpha_k^2\left(2\lambda R^2+ (\lambda+12) \|\nabla f(\bm{x}^{\star})\|^2 + 36R^2 \right) \\
			& \,-(\tfrac{\rho}{2 m \kappa^2} - 2\alpha_k L - \tfrac{4}{\lambda}) \| \bm{x}^k-\Pi_{C}( \bm{x}^k ) \|^2 \\
			&- 2 \alpha_k( f( \Pi_C(\bm{x}^k)) - f^{\star} )+ 48 L^2 \alpha_k^2\|\tilde{\bm{x}}^{l}-\bm{x}^{\star}\|^2.
		\end{split}
	\end{equation*}
\end{lemma}

\begin{proof}
	Since $\bm{x}^{\star}\in X^{\star}\subseteq C$, we have $\bm{x}^{\star}\in  C_0$ and $\bm{x}^{\star}\in  C_{j_k}\subseteq H_k$. 
	Let $\bm{z}^k=\bm{x}^k- \alpha_k\bm{v}^{k}$. 
	It follows that
	\begin{align*}\label{ineq:prop3.2-1}
			&\, \|\bm{x}^{k+1}-\bm{x}^{\star}\|^2
			\leq  \| \Pi_{H_k}\big(\bm{x}^{k}- \alpha_k\bm{v}^{k}\big) - \bm{x}^{\star}\|^2 \\
			\leq &\, \|\bm{x}^k- \alpha_k\bm{v}^{k}-\bm{x}^{\star}\|^2-\|\Pi_{H_k}(\bm{z}^k)-\bm{z}^k\|^2 \\
			= &\, \|\bm{x}^k-\bm{x}^{\star}\|^2 + \alpha_k^2\|\bm{v}^{k}\|^2  \\
			&\, -\|\Pi_{H_k}(\bm{z}^k)-\bm{z}^k\|^2-2 \alpha_k\langle \bm{v}^{k},\bm{x}^{k}-\bm{x}^{\star}\rangle,
	\end{align*}
	where the first inequality follows from Lemma~\ref{lema: Proj-Property}\ref{lema: proj-property-nonexpan} and the fact that $\bm{x}^\star\in C_0$, and the second from the definition of $\bm{z}^k$, Lemma~\ref{lema: Proj-Property}\ref{lema: proj-property-tri-in} and the fact that $\bm{x}^\star\in H_k$.
	Taking conditional expectation on the last inequality and using Lemmas~\ref{lema: SVRG-type-Property} and~\ref{lema: 3.5},
	\begin{equation}
		\label{ineq: prop3.2-2}
		\begin{split}
			& \mathbb{E} \left[\|\bm{x}^{k+1}-\bm{x}^{\star}\|^2\mid \mathcal{F}_k \right]  \\
			\leq \,& \|\bm{x}^k-\bm{x}^{\star}\|^2 - 2 \alpha_k \mathbb{E} \left[ \langle \bm{v}^{k},\bm{x}^{k}-\bm{x}^{\star}\rangle\mid \mathcal{F}_k \right] \\
			\, & + \mathbb{E} \left[ \alpha_k^2 \| \bm{v}^{k}\|^2 - \|\Pi_{H_k} (\bm{z}^k) - \bm{z}^k \|^2 \mid \mathcal{F}_k \right] \\
			\leq \, &     (1+ \alpha_k^2(8L^2+L^2\lambda))\|\bm{x}^k-\bm{x}^{\star}\|^2 \\
			&+(2 \alpha_k L+\tfrac{4}{\lambda})\|\bm{x}^k-\Pi_{C}(\bm{x}^k)\|^2 \\
			&+\alpha_k^2( 2\lambda R^2+(\lambda+4 )\|\nabla f(\bm{x}^{\star})\|^2+12R^2)\\
			&+16L^2 \alpha_k^2\|\tilde{\bm{x}}^{l}-\bm{x}^{\star}\|^2
			-2 \alpha_k (f(\Pi_C(\bm{x}^k))-f^{\star}) \\ 
			&-\mathbb{E} \left[\|\Pi_{H_k}(\bm{z}^k)-\bm{z}^k\|^2\mid\mathcal{F}_k \right].
		\end{split}
	\end{equation}

	We then bound the term  $\mathbb{E}[\|\Pi_{H_k}(\bm{z}^k)-\bm{z}^k\|^2\mid\mathcal{F}_k]$ in last line of \eqref{ineq: prop3.2-2}.
	First, from the triangle inequality and Lemma~\ref{lema: Proj-Property}\ref{lema: proj-property-nonexpan}, we have that
	\begin{equation*}
		\begin{split}
			&\, \|\Pi_{H_k}(\bm{x}^k)-\bm{x}^k\| \\
			\leq & \, \|\bm{z}^k-\bm{x}^k\| 
			+ \|\Pi_{H_k}(\bm{z}^k) -\bm{z}^k\| +  \|\Pi_{H_k}(\bm{x}^k) -\Pi_{H_k}(\bm{z}^k)\|\\
			\leq &\, \|\bm{z}^k-\Pi_{H_k}(\bm{z}^k)\|+2\|\bm{x}^k-\bm{z}^k\|\\
			\leq &\, \|\bm{z}^k-\Pi_{H_k}(\bm{z}^k)\|+2 \alpha_k\|\bm{v}^{k}\|
		\end{split}
	\end{equation*}
	which together with Lemma \ref{lema: SVRG-type-Property} yields that
	\begin{equation*}
		\begin{split}
			&\, \displaystyle \mathbb{E}\left[\|\Pi_{H_k}(\bm{x}^k)-\bm{x}^k\|^2\mid \mathcal{F}_k\right]\\
			\leq&\, \displaystyle 2\mathbb{E}\left[ \| \bm{z}^k-\Pi_{H_k}(\bm{z}^k) \|^2 \mid \mathcal{F}_k \right] + 4 \alpha_k^2 \, \mathbb{E}[ \| \bm{v}^{k} \|^2 \mid \mathcal{F}_k ] \\
			\leq&\, \displaystyle2\mathbb{E}\left[\|\bm{z}^k-\Pi_{H_k}(\bm{z}^k)\|^2\mid\mathcal{F}_k\right]
			+ 32 L^2 \alpha_k^2 \| \bm{x}^k-\bm{x}^{\star} \|^2 \\
			&\,+ 64 L^2 \alpha_k^2 \| \tilde{\bm{x}}^{l}-\bm{x}^{\star}\|^2 + \alpha_k^2( 48 R^2 + 16 \| \nabla f( \bm{x}^{\star}) \|^2).
		\end{split}
	\end{equation*}
	Rearranging the above inequality gives
	\begin{equation}\label{ineq: prop3.2-4}
		\begin{split}
			&-\mathbb{E}\left[\|\bm{z}^k-\Pi_{H_k}(\bm{z}^k)\|^2\mid\mathcal{F}_k \right]\\
			\leq & \, 16 L^2 \alpha_k^2 \|\bm{x}^k-\bm{x}^{\star}\|^2 
			-\frac{1}{2} \mathbb{E} \left [\| \Pi_{H_k}(\bm{x}^k)-\bm{x}^k \|^2 \mid \mathcal{F}_k \right]
			\\
			& \, + 32 L^2\alpha_k^2 \| \tilde{\bm{x}}^{l}-\bm{x}^{\star} \|^2 +  \alpha_k^2( 24 R^2 + 8\|\nabla f(\bm{x}^{\star})\|^2).
		\end{split}
	\end{equation}
	Plugging \eqref{ineq: prop3.2-4} into \eqref{ineq: prop3.2-2} and using Lemma~\ref{delta}, we obtain
	\begin{align*}
		& \mathbb{E}\left[ \|\bm{x}^{k+1} - \bm{x}^{\star} \|^2 \mid \mathcal{F}_k \right] \\
		\leq\, & \left(1 + \alpha_k^2(24 L^2+L^2\lambda)\right)\|\bm{x}^k-\bm{x}^{\star}\|^2  \\
		&\, + \alpha_k^2\left(2\lambda R^2+ (\lambda+12) \|\nabla f(\bm{x}^{\star})\|^2 + 36R^2 \right) \\
		& \,-(\tfrac{\rho}{2 m \kappa^2} - 2\alpha_k L - \tfrac{4}{\lambda}) \| \bm{x}^k-\Pi_{C}( \bm{x}^k ) \|^2 \\
		&\,- 2 \alpha_k ( f( \Pi_C(\bm{x}^k)) - f^{\star} )+ 48 L^2 \alpha_k^2\|\tilde{\bm{x}}^{l}-\bm{x}^{\star}\|^2.
	\end{align*}
	This completes the proof.
\end{proof}

The following lemmas are used in the proof of Theorem~\ref{thm:var_redu_rate}.
\begin{lemma}\label{fea-gap-C-bounded}
Suppose that the conditions of Theorem~\ref{thm:var_redu_rate} hold.
Let $\beta=1-\frac{\rho}{8 m \kappa^2}$.
Then, for any $k\ge 0$, we have 
\begin{align*}
 \mathbb{E}[\|\bm{x}^{k+1}-\Pi_C(\bm{x}^{k+1})\|^2 ]
\le \beta \mathbb{E}[ \|\bm{x}^{k}-\Pi_C (\bm{x}^{k})\|^2 ]
+\mathcal{O}(\alpha_k^2).
\end{align*} 
\end{lemma}

\begin{proof}
Let $\bm{z}^k=\bm{x}^k- \alpha_k\bm{v}^{k}$ and  $\tau= (1-\beta)^{-1}$. 
We have
\begin{equation}\label{xk+1-C-bounded}
\begin{split}
&\, \|\bm{x}^{k+1}-\Pi_C(\bm{x}^{k+1})\|^2 \\
\le &\|\bm{x}^{k+1}-\Pi_C(\bm{x}^k)\|^2 - \|\Pi_{C}(\bm{x}^k) - \Pi_C(\bm{x}^{k+1})\|^2\\
\le &\,\|\Pi_{H_k}(\bm{z}^{k})-\Pi_C(\bm{x}^k)\|^2 - \|\Pi_{C}(\bm{x}^k) - \Pi_C(\bm{x}^{k+1})\|^2\\
\le &\,\|\bm{z}^k-\Pi_C(\bm{x}^k)\|^2-\|\Pi_{H_k}(\bm{z}^k)-\bm{z}^k\|^2\\
\le &\, (1+\tfrac{1}{\tau})\|\bm{x}^k-\Pi_C(\bm{x}^k)\|^2
+(1+\tau)\alpha_k^2\|\bm{v}^k\|^2 \\
&\, -\|\Pi_{H_k}(\bm{z}^k)-\bm{z}^k\|^2,
\end{split}
\end{equation}
where the first inequality follows from Lemma~\ref{lema: Proj-Property}\ref{lema: proj-property-tri-in}, the second from that $\bm{x}^{k+1}, \Pi_{C}(\bm{x}^k)\in C_0$, the third from Lemma~\ref{lema: Proj-Property}\ref{lema: proj-property-tri-in}, and the fourth from the inequality $(a+b)^2 \le (1+ \tfrac{1}{\tau}) a^2 + (1+\tau) b^2$.
We then bound $\mathbb{E}[\|\bm{v}^{k}\|^2\mid\mathcal{F}_{k}]$. 
Using the definition of $\bm{v}^k$, inequality~\eqref{ineq: Lip-cons-f}, Lemma~\ref{lema: Proj-Property}\ref{lema: proj-property-nonexpan} and the boundedness of $C$, we get
\begin{equation}\label{vk-C-bounded}
    \begin{split}
        &\, \mathbb{E}[\|\bm{v}^{k}\|^2\mid\mathcal{F}_{k}] \\
        \le &\, 2L^2\|\bm{x}^k-\tilde{\bm{x}}^l\|^2 + 2\mathbb{E}[\|\nabla f(\tilde{\bm{x}}^{l}))\|^2\mid\mathcal{F}_k ]\\
        \le &\, 4L^2  \|\bm{x}^k-\Pi_{C}(\bm{x}^k)\|^2+ 4L^2\|\Pi_{C}(\bm{x}^k)-\tilde{\bm{x}}^l\|^2  \\
        & \, + 4 \|\nabla f( \tilde{\bm{x}}^l) - \nabla f(\Pi_C(\tilde{\bm{x}}^l))\|^2 + 4 \|\nabla f(\Pi_C(\tilde{\bm{x}}^l))\|^2\\
        \le &\, 4L^2\|\bm{x}^k-\Pi_{C}(\bm{x}^k)\|^2+12L^2\|\Pi_{C}(\tilde{\bm{x}}^l)-\tilde{\bm{x}}^l\|^2
        +\mathcal{O}(1).
    \end{split}
\end{equation} 
We next bound
$-\mathbb{E}[\|\bm{z}^k-\Pi_{H_k}(\bm{z}^k)\|^2\mid\mathcal{F}_k]$.
Using Lemma~\ref{lema: Proj-Property}\ref{lema: proj-property-nonexpan},
\begin{align*}
&\, \|\Pi_{H_k}(\bm{x}^k)-\bm{x}^k\| \\
\le &\, \|\Pi_{H_k}(\bm{x}^k) -\Pi_{H_k}(\bm{z}^k)\|+ \|\Pi_{H_k}(\bm{z}^k) -\bm{z}^k\|+\|\bm{z}^k-\bm{x}^k\| \\
\le &\, \|\bm{z}^k-\Pi_{H_k}(\bm{z}^k)\|+2\|\bm{x}^k-\bm{z}^k\|\\
= &\, \|\bm{z}^k-\Pi_{H_k}(\bm{z}^k)\|+2 \alpha_k\|\bm{v}^{k}\|.
\end{align*}
Taking conditional expectation and using inequality~\eqref{vk-C-bounded},
\begin{equation*}
\begin{split}
    &\, \displaystyle \mathbb{E}[\|\Pi_{H_k}(\bm{x}^k)-\bm{x}^k\|^2\mid \mathcal{F}_k]\\
    \le&\, \displaystyle2\mathbb{E}[\|\bm{z}^k-\Pi_{H_k}(\bm{z}^k)\|^2\mid\mathcal{F}_k]  + 16^2\alpha_k^2 L^2 \|\bm{x}^k-\Pi_{C}(\bm{x}^k)\|^2 \\
    &\, +48 \alpha_k^2 L^2 \|\Pi_{C}(\tilde{\bm{x}}^l)-\tilde{\bm{x}}^l\|^2
        +\mathcal{O}(\alpha_k^2).
\end{split}
\end{equation*}
Rearranging the above inequality and using Lemma~\ref{delta} yields
\begin{equation}\label{zk-C-bounded}
\begin{split}
    &\, \displaystyle -\mathbb{E}[\|\bm{z}^k-\Pi_{H_k}(\bm{z}^k)\|^2\mid\mathcal{F}_k]\\
    \le& \, \left(-\tfrac{\rho}{2m\kappa^2}+8 L^2\alpha_k^2\right) \|\bm{x}^k-\Pi_{C}(\bm{x}^k)\|^2\\
    &\, + 24L^2 \alpha_k^2\|\Pi_{C}(\tilde{\bm{x}}^l)-\tilde{\bm{x}}^l\|^2 
        +\mathcal{O}(\alpha_k^2).
\end{split}
\end{equation}
Using \eqref{xk+1-C-bounded}, \eqref{vk-C-bounded} and~\eqref{zk-C-bounded} and the definitions of $\tau$ and $\alpha_k$,
\begin{equation}\label{dis-xk+1-C-bounded}
\begin{split}
&\mathbb{E}[\|\bm{x}^{k+1}-\Pi_C(\bm{x}^{k+1})\|^2\mid \mathcal{F}_k]\\
\le  & \, (1 -\tfrac{\rho}{4m\kappa^2})\|\bm{x}^k-\Pi_C(\bm{x}^k)\|^2   +\mathcal{O}(\alpha_k^2) \\
&\, +12L^2 \alpha_k^2(3+\tau) \|\Pi_{C}(\tilde{\bm{x}}^l)-\tilde{\bm{x}}^l\|^2. 
   \end{split}
\end{equation}
Using the tower property of conditional expectation~\cite[Theorem~2.3.2(iii)]{Shreve2004StochasticCF}, (\ref{dis-xk+1-C-bounded}) and the boundedness of $\alpha_k$, we have
\begin{equation*}
\begin{split}
&\, \mathbb{E}[\|\tilde{\bm{x}}^{l+1}-\Pi_C(\tilde{\bm{x}}^{l+1})\|^2\mid \mathcal{F}_{lr}]\\
\le \; &(1-\tfrac{\rho}{4m\kappa^2})\|\tilde{\bm{x}}^l-\Pi_{C}(\tilde{\bm{x}}^l)\|^2+\mathcal{O}( \alpha_{lr}^2)\\
&\, + 12rL^2\alpha_k^2(3+\tau)\|\Pi_{C}(\tilde{\bm{x}}^l)-\tilde{\bm{x}}^l\|^2
\\
\le \, &(1-\tfrac{\rho}{8m\kappa^2})\|\tilde{\bm{x}}^l-\Pi_{C}(\tilde{\bm{x}}^l)\|^2
+\mathcal{O}(1)\\
\le \, & \beta^{l+1}\|\bm{x}^0-\Pi_{C}(\bm{x}^0)\|^2 + \mathcal{O}(1),
   \end{split}
\end{equation*}
where the last inequality follows from inductively using the second last one.
Therefore, the sequence $\{\mathbb{E}[\|\tilde{\bm{x}}^l-\Pi_{C}(\tilde{\bm{x}}^l)\|^2\mid \mathcal{F}_{lr}]\}_l$ is bounded. By~\eqref{dis-xk+1-C-bounded}, we get the desired inequality.
\end{proof}

\begin{lemma}\label{lem:x^k_bound}
Suppose that the conditions of Theorem~\ref{thm:var_redu_rate} hold. Then, $\{\mathbb{E}[\|\bm{x}^k-\bm{x}^{\star}\|]\}_k$ and $\{\mathbb{E}[\|\bm{x}^{k}\|]\}_k$ are bounded sequences.
\end{lemma}
\begin{proof}
For any $\bm{x}^\star\in X^\star$, applying Lemma~\ref{fea-gap-C-bounded} inductively yields
\begin{equation*}
\begin{split}
&\, \mathbb{E}[\|\bm{x}^{k+1}-\bm{x}^{\star}\|^2 ]\\
\le &\, 2\mathbb{E}[\|\bm{x}^{k+1}-\Pi_C(\bm{x}^{k+1})\|^2 ]
+ 2\mathbb{E}[\|\Pi_C(\bm{x}^{k+1})-\bm{x}^{\star}\|^2 ]\\
\le  &\, 2\beta^{k+1}\|\bm{x}^0-\Pi_{C}(\bm{x}^0)\|^2 + \mathcal{O}(1) + 2 \mathbb{E}[\|\Pi_C(\bm{x}^{k+1})-\bm{x}^{\star}\|^2 ].
\end{split}
\end{equation*}
Since $C$ is bounded, the sequences
$\{\mathbb{E}[\|\bm{x}^k-\bm{x}^{\star}\|^2 ]\}$ and $\{\mathbb{E}[\|\bm{x}^k\|^2 ]\}$ are bounded. The proof is completed.
\end{proof}

\begin{lemma}\label{x^k+1-tilde_x^l}
Suppose that the conditions of Theorem~\ref{thm:var_redu_rate} hold. Then, for $k\in \{lr,lr+1,\cdots,lr+r-1\}$, we have
\begin{equation*}
\mathbb{E}[\|\bm{x}^{k+1}-\tilde{\bm{x}}^l\|^2 ]
\le \mathcal{O}(\alpha_k^2)
+2r \sum_{t=lr}^{t=k}\mathbb{E}[\|\bm{x}^t-\Pi_C(\bm{x}^t)\|^2 ].
\end{equation*}
\end{lemma}

\begin{proof}
We first prove that 
\begin{equation}\label{ineq:proof:9}
\mathbb{E}[\|\bm{x}^{k+1}-\bm{x}^k\|^2]
\le \mathcal{O}(\alpha_k^2) +2\mathbb{E}[\|\bm{x}^k-\Pi_C(\bm{x}^k)\|^2].
\end{equation}
Since $C_0 = \mathbb{R}^d$,
\begin{equation*}
\bm{x}^{k+1}=
\begin{cases}
  \bm{x}^{k}-\alpha_k \bm{v}^{k}-\tfrac{ \left( \phi_{j_k}(\bm{x}^k)-\alpha_k\langle \bm{\xi}^k, \bm{v}^{k}\rangle \right)_+ }{\|\bm{\xi}^k\|^2}\bm{\xi}^k& \text{if } \bm{\xi}^k \neq \bm{0},\\
\bm{x}^{k}-\alpha_k \bm{v}^{k} &\text{if } \bm{\xi}^k = \bm{0}.
\end{cases}
\end{equation*}
If $\bm{\xi}^k=\bm{0}$, it follows from \eqref{vk-C-bounded} and Lemma~\ref{lem:x^k_bound} that
\begin{align*}
\mathbb{E}[\|\bm{x}^{k+1}-\bm{x}^k\|^2 ]
= \alpha_k^2\,\mathbb{E}[\|\bm{v}^k\|^2 ]
\le \mathcal{O}(\alpha_k^2), 
\end{align*}
which implies \eqref{ineq:proof:9}. If $\bm{\xi}^k\neq \bm{0}$, then
\begin{equation*}\label{x{k+1}-x{k}}
\begin{split}
&\, \|\bm{x}^{k+1}-\bm{x}^k\| 
= \|\alpha_k \bm{v}^k + \tfrac{ \left( \phi_{j_k}(\bm{x}^k)-\alpha_k\langle \bm{\xi}^k, \bm{v}^{k}\rangle \right)_+ }{\|\bm{\xi}^k\|^2}\bm{\xi}^k \| \\
\le &\, \alpha_k \|\bm{v}^k\| + \tfrac{(\phi_{j_k}(\bm{x}^k))_{+}}{\|\bm{\xi}^k\|}+ \tfrac{ (\alpha_k\langle \bm{\xi}^k, \bm{v}^{k}\rangle )_+ }{\|\bm{\xi}^k\|} \\
\le &\, 2\alpha_k \|\bm{v}^k\| + \|\bm{x}^k - \Pi_{H_k}(\bm{x}^k)\| \\
\le&\,  2\alpha_k \|\bm{v}^k\|  +  \|\bm{x}^k - \Pi_{C}(\bm{x}^k)\|,
\end{split}
\end{equation*}
where the second inequality follows from Lemma~\ref{halfspace}, and
the last from $\Pi_{C}(\bm{x}^k) \in C\subseteq H_k$.
By \eqref{vk-C-bounded} and Lemma~\ref{lem:x^k_bound}, \eqref{ineq:proof:9} follows. For $k\in \{lr,lr+1,\cdots,lr+r-1\}$, \eqref{ineq:proof:9} implies
\begin{equation*}
\begin{split}
&\, \mathbb{E}[\|\bm{x}^{k+1}-\tilde{\bm{x}}^l\|^2 ]
\le r\sum_{t=lr}^{t=k}\mathbb{E}[\|\bm{x}^{t+1}-\bm{x}^t\|^2 ]\\
\le&\, \mathcal{O}(\alpha_k^2)
+2r \sum_{t=lr}^{t=k}\mathbb{E}[\|\bm{x}^t-\Pi_C(\bm{x}^t)\|^2 ].
\end{split}
\end{equation*}
The proof is completed.
\end{proof}

\subsection{Proof of Proposition~\ref{prop:LR-holds}}\label{appendix-A}

The proof  of case~\ref{LR-sufficient-poly} can be founded in  \cite{hoffman2003approximate}.
We therefore prove only case~\ref{LR-sufficient-bounded}.
If $\bm{x}\in C$, then \eqref{LR-condition} trivially holds.
We thus assume that $\bm{x} \notin C$. 
Let $I(\bm{x})$ be the set defined by
\begin{equation}
	\label{def:I(x)}
	I(\bm{x})=\{j\in[m] \mid  \phi_j(\bm{x})>0\}.
\end{equation}
Since $\bm{x} \in C_0$, we have that $\bm{x}\notin \cap_{j\in[m]}C_j$ and hence the index set $I(\bm{x})$ is non-empty.
By~\cite[Theorem~9]{bauschke1999strong},
there exists  a constant $\gamma>0$ such that
\begin{equation}\label{ineq: result-lema-LR-Bounded}
	\dist(\bm{x},C) \leq \gamma \phi_{j'}(\bm{x}),
\end{equation}
where $j'\ \in \arg\max_{j\in[m]}  \phi_{j}(\bm{x})  $ and $\gamma$ does not depend on $\bm{x}$.
Clearly, 
$j'\in I(\bm{x})$  and hence $\phi_{j'}(\bm{x})=\max_{j\in I(\bm{x})}\phi_{j}(\bm{x})$.
By supposition,
there exists a constant $\eta>0$ such that for any $j\in [m]$, $\bm{x}\in C_0$ and $\bm{\xi}_{j}\in \partial \phi_j (\bm{x})$, we have $\|\bm{\xi}_{j}\|\leq \eta$.
Fix an arbitrary $\bm{\xi}_{j^{'}}\in  \partial \phi_{j'} (\bm{x})$. If $\bm{\xi}_{j^{'}} \neq \bm{0}$,
\begin{equation*}
	\begin{split}
		\phi_{j'}(\bm{x}) \le\, & -\langle \bm{\xi}_{j'}, \Pi_{H(\phi_{j'};\, \bm{x};\,\bm{\xi}_{j'}))}(\bm{x})-\bm{x}\rangle \\
		\leq \, & \eta \dist(\bm{x},\,H(\phi_{j'};\, \bm{x};\,\bm{\xi}_{j'})),
	\end{split}
\end{equation*}
where the first inequality follows from the definition of $H(\phi_{j'};\, \bm{x};\,\bm{\xi}_{j'})$ and the fact that~$\Pi_{H(\phi_{j'};\, \bm{x};\,\bm{\xi}_{j'}))}(\bm{x}) \in H(\phi_{j'};\, \bm{x};\,\bm{\xi}_{j'})$, and the second inequality follows from the Cauchy-Schwarz inequality and the uniform bound of the subdifferential $ \phi_j (\bm{x})$.
If $\bm{\xi}_{j^{'}} = \bm{0}$, by convexity, $\bm{x}$ is a minimizer of $\phi_{j'}$ and 
\[ \phi_{j'}(\bm{x}) = \min_{\bm{y}\in\mathbb{R}^d} \phi_{j'}(\bm{y}) \le 0,\]
where the inequality follows from the non-emptiness of $\{ \bm{y}\in \mathbb{R}^d \mid \phi_j(\bm{y}) < 0,\ j \in [m] \} $.
However, by the definition of $j'$, $\phi_{j'}(\bm{x}) > 0$. Hence, it is impossible to have $\bm{\xi}_{j^{'}} = \bm{0}$.
Therefore, 
\[\phi_{j'}(\bm{x})  \leq   \eta \min_{\bm{\xi}_{j'}\in   \partial \phi_{j'} (\bm{x}) }\ \dist(\bm{x},H(\phi_{j'};\, \bm{x};\,\bm{\xi}_{j'})),\]
which, together with \eqref{ineq: result-lema-LR-Bounded}, implies that
\begin{equation*}
	\begin{split}
		&\, \dist(\bm{x},C)  \\  
		\leq       & \,  \gamma \eta \min_{\bm{\xi}_{j'}\in   \partial \phi_{j'} (\bm{x}) }\ \dist(\bm{x},H(\phi_{j'};\, \bm{x};\,\bm{\xi}_{j'})) \\
		\le       &\,  \gamma \eta \max_{j\in [m]} \min_{\bm{\xi}_{j}\in \partial \phi_j (\bm{x})}
		\dist(\bm{x},\, H(\phi_{j};\, \bm{x};\,\bm{\xi}_{j})). 
	\end{split}
\end{equation*}
Noting that both constants $\gamma$ and $\eta$ are independent of $\bm{x}$, this completes the proof.

\subsection{Proof of Theorem~\ref{convergence1}}\label{subsec:converge}
We first prove that the sub-sequence $\{\tilde{\bm{x}}^l\}$ converges.
Since $\mu_l\to 0$ as $l\to \infty$, 
there exists $l_0\geq 0$ such that $2 \mu_l L \leq \frac{\rho}{8m \kappa^2} $ for any $l\geq l_0$.
Take $\lambda = 32 m \kappa^2\rho ^{-1}$. Then, for any for any $l\geq l_0$,
\[ \frac{\rho}{2m \kappa^2}-2 \mu_l L-\frac{4}{\lambda}\geq \frac{\rho}{4m\kappa^2}>0. \]
Fix any optimal solution $\bm{x}^{\star} \in X^\star$. It follows from Lemma~\ref{lemma: one-iteration} and the definition of $\alpha_k$ that for all $l\ge l_0$,
\begin{align}
		&\mathbb{E}\left[\|\bm{x}^{k+1}-\bm{x}^{\star}\|^2\mid\mathcal{F}_{k}\right] \notag\\
		\leq & \, ( 1+\mathcal{O}( \alpha_k^2)) \,\| \bm{x}^k - \bm{x}^{\star} \|^2 +    \mathcal{O}( \alpha_k^2)       + \mathcal{O}( \alpha_k^2)   \| \tilde{\bm{x}}^{l} - \bm{x}^{\star} \|^2  \notag  \\
		& \,- \frac{\rho}{4m\kappa^2} \| \bm{x}^k-\Pi_{C}(\bm{x}^k) \|^2 
		- 2 \mu_l ( f(\Pi_C(\bm{x}^k) ) - f^{\star} )  \label{ineq: theorem3.1-1} \\ 
		\leq &\, (1 + \mathcal{O} ( \mu_l^2)) \, \|\bm{x}^k - \bm{x}^{\star}\|^2 + \mathcal{O} ( \mu_l^2 )\|\tilde{\bm{x}}^{l} - \bm{x}^{\star}\|^2 + \mathcal{O}( \mu_l^2).\notag
\end{align}
Using the tower property of conditional expectation~\cite[Theorem~2.3.2(iii)]{Shreve2004StochasticCF} and inequality~\eqref{ineq: theorem3.1-1},
\begin{align}
	&\, \mathbb{E} \left[\| \bm{x}^{lr+r} - \bm{x}^{\star}  \|^2\mid \mathcal{F}_{lr} \right] \notag \\
	= &\, \mathbb{E}\left[ \mathbb{E} \left[\| \bm{x}^{lr+r-1 + 1} - \bm{x}^{\star}  \|^2\mid \mathcal{F}_{lr+r - 1} \right] \Big\vert \mathcal{F}_{lr} \right]  \notag \\
	\le &\, (1 + \mathcal{O} ( \mu_l^2) ) \mathbb{E}\left[ \|\bm{x}^{lr + r - 1} - \bm{x}^{\star}\|^2  \Big\vert \mathcal{F}_{lr} \right]  \notag \\
	&\, + \mathcal{O} ( \mu_l^2 )\|\tilde{\bm{x}}^{l} - \bm{x}^{\star}\|^2 + \mathcal{O}( \mu_l^2) \notag \\
	&\vdots \label{ineq:proof:6}\\
	\leq & \, (1+\mathcal{O}(\mu_l ^2))^{r-1} \mathbb{E} [\|\bm{x}^{lr + 1} - \bm{x}^{\star}\|^2\mid \mathcal{F}_{lr} ]  \notag \\
	&\, +  \left(\mathcal{O}( \mu_l^2)\| \tilde{\bm{x}}^{l} - \bm{x}^{\star}\|^2   \right)   \sum_{s=0}^{r-2} (1 + \mathcal{O}(\mu_l ^2))^{s}  \notag \\
	&\, +  \mathcal{O}( \mu_l^2)   \sum_{s=0}^{r-2} (1 + \mathcal{O}(\mu_l ^2))^{s}. \notag 
\end{align}
Similarly, we get
\begin{equation}
	\begin{split}
		&\, \mathbb{E} \left[\|\bm{x}^{lr + 1} - \bm{x}^{\star} \|^2 \mid \mathcal{F}_{lr} \right]   \\
		\leq & \, ( 1+\mathcal{O}( \mu_l^2)) \,\| \tilde{\bm{x}}^{l} - \bm{x}^{\star} \|^2 - \frac{\rho}{4m\kappa^2} \| \tilde{\bm{x}}^{l} -\Pi_{C}(\tilde{\bm{x}}^{l} ) \|^2     \\
		& \,  + \mathcal{O} ( \mu_l^2 )  - 2 \mu_l ( f(\Pi_C(\tilde{\bm{x}}^{l}) ) - f^{\star} ).\label{ineq:proof:5}
	 \end{split}
\end{equation}
Since $\mu_l < 1$ for sufficiently large $l$, $(1+\mathcal{O}(\mu_l ^2))^{s} \le 1 + \mathcal{O}(\mu_l^2)$ for any $l$ and $s\in [r]$.
We thus have
\begin{align}
	&\, \mathbb{E} \left[\| \tilde{\bm{x}}^{(l+1)} - \bm{x}^{\star}  \|^2\mid \mathcal{F}_{lr} \right]  \notag \\
	\leq\, & (1+\mathcal{O}(\mu_l ^2))^r \mathbb{E} \left[\|\bm{x}^{lr + 1}  - \bm{x}^{\star}\|^2\mid \mathcal{F}_{lr} \right]  \notag \\
	& +  \left(\mathcal{O}( \mu_l^2)\| \tilde{\bm{x}}^{l} - \bm{x}^{\star}\|^2 +  \mathcal{O}( \mu_l^2)  \right)   \sum_{s=0}^{r-1} (1 + \mathcal{O}(\mu_l ^2))^{s}   \label{ineq: theorem3.1-7} \\
	&\, - \frac{\rho}{4m\kappa^2} \| \tilde{\bm{x}}^{l} -\Pi_{C}(\tilde{\bm{x}}^{l} ) \|^2  - 2 \mu_l ( f(\Pi_C(\tilde{\bm{x}}^{l}) ) - f^{\star} )  \notag \\
	\leq&\,(1+\mathcal{O}(\mu_l^2))\,\|\tilde{\bm{x}}^{l}-\bm{x}^{\star}\|^2 +\mathcal{O}( \mu_l^2 ) \notag \\
	&\, - \frac{\rho}{4m\kappa^2}  \|\tilde{\bm{x}}^{l}-\Pi_{C}(\tilde{\bm{x}}^{l})\|^2 , \notag 
\end{align}
where the first inequality follows by substituting~\eqref{ineq:proof:5} into~\eqref{ineq:proof:6}. Note that now the constants hidden in the big-O notation could possibly depend on $r$, but they are independent of $\mu_l$ or $l$. Applying Lemma~\ref{lema: supermartingale} to the recursion~\eqref{ineq: theorem3.1-7}, we have that the sequence $\{\|\tilde{\bm{x}}^{l}-\bm{x}^{\star}\|^2\}$ converges almost surely, that 
\begin{equation}\label{ineq:proof:7}
	\sum_{l=0}^{\infty} \mu_l [f( \Pi_{C}( \tilde{\bm{x}}^l ) ) - f^{\star}] < \infty\quad a.s.,
\end{equation}
and that
\begin{equation}\label{ineq:proof:8}
	\sum_{l=0}^{\infty}   \|\tilde{\bm{x}}^{l} - \Pi_{C}( \tilde{\bm{x}}^{l} ) \|^2 < \infty\quad a.s.
\end{equation}
By inequality~\eqref{ineq:proof:7} and the fact that $\sum_{l=0}^{\infty} \mu_l=\infty$,
\begin{equation}\label{ineq: theorem3.1-8}
	\liminf_{l\to \infty} f(\Pi_{C}(\tilde{\bm{x}}^l)) =f^{\star} \quad a.s.
\end{equation}
Also, inequality~\eqref{ineq:proof:8} implies that
\begin{equation}\label{ineq: theorem3.1-9}
	\lim_{l\rightarrow \infty}\|\Pi_{C}(\tilde{\bm{x}}^l)- \tilde{\bm{x}}^l\|=0 \quad a.s.
\end{equation}
Since the sequence $\{\|\tilde{\bm{x}}^l - \bm{x}^{\star} \|\}$ converges almost surely, the sequence $\{\tilde{\bm{x}}^l\}$ is bounded and has an accumulation point $\tilde{\bm{x}}^{\star}$ almost surely. 
Therefore, there exists a sub-sequence $\{\tilde{\bm{x}}^{l_t}\}$ such that $ \tilde{\bm{x}}^{l_t} \to \tilde{\bm{x}}^{\star}$ as $t\to \infty$. 
By relation~\eqref{ineq: theorem3.1-9} and continuity of $\Pi_{C}(\cdot)$, the sequence $ \Pi_{C}( \tilde{\bm{x}}^{l_t} ) $ converges almost surely to $\Pi_{C}( \tilde{\bm{x}}^{\star} ) = \tilde{\bm{x}}^{\star}\in C$. 
It follows from \eqref{ineq: theorem3.1-8} and the continuity of $f$ that $f(\tilde{\bm{x}}^{\star}) = f^{\star}$. 
Hence,  $\tilde{\bm{x}}^{\star}\in X^{\star}$. Since $\| \tilde{\bm{x}}^{l} - \bm{x}^{\star} \|$ converges almost surely for every $\bm{x}^{\star}\in X^{\star}$, we have that $\| \tilde{\bm{x}}^{l} - \tilde{\bm{x}}^{\star} \|$ converges almost surely. 
Since $\| \tilde{\bm{x}}^{l_t} - \tilde{\bm{x}}^{\star} \| \to 0$ as $t \to \infty$ almost surely, we have that $\| \tilde{\bm{x}}^{l} - \tilde{\bm{x}}^{\star} \| \to 0$ as $l\to \infty$ almost surely. Thus, almost surely, we have $\lim_{l\to \infty}\tilde{\bm{x}}^{l} =\tilde{\bm{x}}^{\star} $. 
To prove the convergence in $\{\bm{x}^k\}$, by the boundedness of the sequence~$\{ \|\tilde{\bm{x}}^{l}-\tilde{\bm{x}}^{\star}\|^2 \}$ and Lemma~\ref{lemma: one-iteration},
\begin{equation*} 
	\begin{split}
		&\, \mathbb{E} \left[\|\bm{x}^{k+1}-\tilde{\bm{x}}^{\star}\|^2\mid\mathcal{F}_{k} \right]  \\
		\leq\,& (1+ \mathcal{O}(\alpha_k^2)) \| \bm{x}^k - \tilde{\bm{x}}^{\star} \|^2 + \mathcal{O}(\alpha_k^2)  \\
		& \, - \frac{\rho}{4m\kappa^2} \|\bm{x}^k-\Pi_{C}(\bm{x}^k)\|^2  - 2 \alpha_k(f(\Pi_C(\bm{x}^k)) - f^{\star}),
	\end{split}
\end{equation*}
which, together with 
Lemma~\ref{lema: supermartingale} and the fact that $\sum_{k=0}^{\infty}\alpha_k^2< \infty$, implies that the sequence $\{\|\bm{x}^{k}-\tilde{\bm{x}}^{\star}\|^2\}$ converges almost surely.
Since the sub-sequence $\{\|\tilde{\bm{x}}^{l}-\tilde{\bm{x}}^{\star}\|^2\}$ converges almost surely to $0$, we have that $\{\|\bm{x}^{k}-\tilde{\bm{x}}^{\star}\|^2\}$ converges almost surely to $0$ as well, which shows that $\lim_{k \to \infty} \bm{x}^k = \tilde{\bm{x}}^{\star}$.
This completes the proof.

\subsection{Proof of Theorem~\ref{rate_pro1}}\label{subsec:converge_rate}
We first prove the convergence rate of the feasibility gap.
Fix an arbitrary optimal solution $\bm{x}^{\star} \in X^\star$.
By using Lemma~\ref{lemma: one-iteration} with $\lambda=16m\kappa^2\rho^{-1}$ and the definition of $\alpha_k$,  we have that for all $k\geq0$, 
\[ \frac{\rho}{2m \kappa^2}-2\alpha_k L-\frac{4}{\lambda} \geq \frac{\rho}{8m\kappa^2},\] 
and hence that
\begin{align}
	& \mathbb{E}\left[\|\bm{x}^{k+1}-\bm{x}^{\star}\|^2\mid\mathcal{F}_{k}
	\right] \notag \\ 
	\leq \, &( 1 + \mathcal{O}( \alpha_k^2)) \|\bm{x}^k - \bm{x}^{\star}\|^2 + \mathcal{O}( \alpha_k^2) \|\tilde{\bm{x}}^{l} -\bm{x}^{\star}\|^2    \notag  \\
	& \, + \mathcal{O}( \alpha_k^2 ) 	- \frac{\rho}{8m\kappa^2} \|\bm{x}^k - \Pi_{C} (\bm{x}^k) \|^2  \label{ineq: theorem3.2-1}\\
	& \,  - 2 \alpha_k (f(\Pi_C(\bm{x}^k)) - f^{\star} ) \notag \\
	\leq \, & ( 1 + \mathcal{O}( \alpha_k^2)) \|\bm{x}^k - \bm{x}^{\star}\|^2  + \mathcal{O}( \alpha_k^2) \|\tilde{\bm{x}}^{l} - \bm{x}^{\star} \|^2 + \mathcal{O}( \alpha_k^2 ). \notag 
\end{align}
Since $C_0$ is compact, the sequences $\{\|\bm{x}^k - \bm{x}^{\star}\|^2 \}_k$ and  $\{\|\tilde{\bm{x}}^l - \bm{x}^{\star}\|^2 \}_l $ are bounded. By inequality~\eqref{ineq: theorem3.2-1}, for any $k\ge 0$,
\begin{align}
	& \frac{\rho}{8m\kappa^2} \mathbb{E} \left[\|\bm{x}^k - \Pi_{C}(\bm{x}^k)\|^2 \right] \notag\\ 
	\leq\,    &  \mathbb{E} \left[\|\bm{x}^k - \bm{x}^{\star}\|^2 \right]  + \mathcal{O}( \alpha_k^2 ) \notag\\
	&\, - \mathbb{E} \left[ \mathbb{E} \left[\|\bm{x}^{k+1} - \bm{x}^{\star} \|^2 \mid\mathcal{F}_{k} \right] \right] \label{ineq: theorem3.2-2}\\
	= \,  & \mathbb{E} \left[\|\bm{x}^k - \bm{x}^{\star}\|^2 \right]
	-   \mathbb{E} \left[\|\bm{x}^{k+1} - \bm{x}^{\star} \|^2 \right]   + \mathcal{O}( \alpha_k^2 ). \notag
\end{align}
Summing the last inequality over $k$, we get
\begin{equation*}
	\begin{split}
		&  \frac{\rho}{8m\kappa^2} \sum_{k=0}^{K-1} \mathbb{E} \left[ \|\bm{x}^k - \Pi_{C}(\bm{x}^k) \|^2 \right] \\
		\leq & \, \| \bm{x}^0 - \bm{x}^{\star} \|^2 - \mathbb{E} \left[ \|\bm{x}^{K} - \bm{x}^{\star}\|^2 \right] + \mathcal{O}(1) \sum_{k=0}^{K-1} \alpha_k^2 \\
		\leq & \, \mathcal{O}(1)+\mathcal{O}(1) \sum_{k=1}^{K} \frac{1}{k} \leq \mathcal{O}(1) \log(K). 
	\end{split}
\end{equation*}
By the convexity of $\dist^2(\cdot,C)$, it follows that
\begin{equation*}
	\begin{split}
		\mathbb{E} \left[\dist^2(\bar{\bm{x}}^{K}, C) \right] 
		\leq \, & \frac{1}{K} \sum_{k=0}^{K-1} \mathbb{E} \left[\|\bm{x}^k - \Pi_{C}(\bm{x}^k)\|^2 \right]\\
		\leq \, & \mathcal{O} \left( \frac{\log (K)}{K} \right).  
	\end{split}
\end{equation*}
Next, we prove the convergence rate of the optimality gap. By the definition of $\bm{v}^{k}$, we have
\begin{equation*}
	\mathbb{E} \left[  \langle \bm{v}^{k}, \bm{x}^k - \bm{x}^{\star} \rangle \mid \mathcal{F}_{k} \right]
	\geq f(\bm{x}^k) - f^{\star} ,
\end{equation*}
which, together with that $\bm{x}^\star\in H_k, C_0$ and Lemma~\ref{lema: SVRG-type-Property}, implies
\begin{align*}
	&\, \mathbb{E} \left[ \|\bm{x}^{k+1} - \bm{x}^{\star} \|^2 \mid \mathcal{F}_k \right] \le \mathbb{E} \left[ \|\bm{x}^k - \alpha_k \bm{v}^k - \bm{x}^{\star} \|^2 \mid \mathcal{F}_k \right]\\
	= &\, \mathbb{E} \left[ \|\bm{x} - \bm{x}^\star\|^2 + \alpha_k^2 \|\bm{v}^k\|^2 - 2\alpha_k \langle \bm{x}^k - \bm{x}^\star , \bm{v}^k \rangle\mid \mathcal{F}_k \right]\\
	\leq &\,  ( 1 + 8L^2 \alpha_k^2 )\|\bm{x}^k - \bm{x}^{\star} \|^2 + 16L^2 \alpha_k^2 \| \tilde{\bm{x}}^{l} - \bm{x}^{\star} \|^2\\
	&\,- 2 \alpha_k \left( f(\bm{x}^k) - f^{\star} \right) + 4 \| \nabla f(\bm{x}^{\star}) \|^2 \alpha_k^2   + 12 R^2 \alpha_k^2  \\
	\le &\, \|\bm{x}^k - \bm{x}^{\star} \|^2 - 2 \alpha_k \left( f(\bm{x}^k) - f^{\star} \right) + \mathcal{O}(1) \alpha_k^2,
\end{align*}
where the last inequality follows from the boundedness of $C_0$ and the fact that~$\tilde{\bm{x}}^{l}, \bm{x}^{k}, \bm{x}^{\star}\in C_0$.
Taking expectation on both sides, we obtain
\begin{equation*}
	\begin{split}
		&\,  \mathbb{E} \left[ f(\bm{x}^k) - f^{\star} \right] \\
		\leq &\,  \frac{1}{2\alpha_k} \left(\mathbb{E} \left[ \|\bm{x}^{k} - \bm{x}^{\star} \|^2 \right] - \mathbb{E} \left[ \|\bm{x}^{k+1} - \bm{x}^{\star} \|^2 \right] \right)  + \mathcal{O}(\alpha_k).
	\end{split}
\end{equation*}
Summing of the last inequality over $k$, we get
\begin{align}
	&\, \sum_{k=0}^{K-1}\mathbb{E} \left[ f(\bm{x}^k) - f^{\star} \right] \notag \\
	\le &\, \frac{1}{2\tilde{\alpha} } \mathbb{E}\left[ \| \bm{x}^0 - \bm{x}^{\star} \|^2 \right] - \frac{\sqrt{K} }{2\tilde{\alpha} } \mathbb{E}\left[ \| \bm{x}^K - \bm{x}^{\star} \|^2 \right] \notag  \\
	&\, + \frac{1}{2\tilde{\alpha}}\sum_{k\in [K-1]}\left( \sqrt{k+1} - \sqrt{k} \right) \mathbb{E} \left[ \|\bm{x}^{k} - \bm{x}^{\star} \|^2 \right] \label{ineq: theorem3.2-3}  \\
	&\, + \mathcal{O}\left( \sum_{k\in [K]} \frac{1}{\sqrt{k}} \right) \notag\\
	\leq &\, \mathcal{O}(1)+ \mathcal{O}(1)\sum_{k\in [K-1]} \left( \sqrt{k+1} - \sqrt{k} \right)   + \mathcal{O}\left(\sum_{k\in [K]} \frac{1}{\sqrt{k}} \right) \notag  \\
	\leq &\, \mathcal{O}(1) + \mathcal{O}(\sqrt{K - 1}) + \mathcal{O}(\sqrt{K}) = \mathcal{O}(\sqrt{K} ), \notag 
\end{align}
where the second inequality follows from the boundedness of $C_0$, and the last from the fact that 
\[  \sum_{k\in [K]} \frac{1}{\sqrt{k}} \le \int_1^{K+1}  \frac{dt}{\sqrt{t}}  = 2(\sqrt{K+1} -1 ) .\]
By the convexity of $f$ and \eqref{ineq: theorem3.2-3}, we have 
\begin{equation*}
	\mathbb{E} [f(\bar{\bm{x}}^{K})]\leq \frac{1}{K} \sum_{k=0}^{K-1} \mathbb{E}[f(\bm{x}^k)] \leq f^{\star}+\mathcal{O} \left( \frac{1}{\sqrt{K}} \right).
\end{equation*}
This completes the proof.

\subsection{Proof of Theorem~\ref{rate_convex_constant_stepsize}}\label{subsec:converge_rate_constant_step}
Similarly to~\eqref{ineq: theorem3.2-2}, we get for any $k\ge 0$,
\begin{equation*} 
	\begin{split}
		&	\frac{\rho}{8m\kappa^2} \mathbb{E} \left[\|\bm{x}^k - \Pi_{C}(\bm{x}^k)\|^2 \right]  \\
		\leq\, & \mathbb{E} \left[\|\bm{x}^k - \bm{x}^{\star}\|^2 \right]
		+ \mathcal{O}( K^{-1} )\\
		& \,  - \mathbb{E} \left[ \mathbb{E} \left[\|\bm{x}^{k+1} - \bm{x}^{\star} \|^2 \mid\mathcal{F}_{k} \right] \right]\\ 
		= \,	&  \mathbb{E} \left[\|\bm{x}^k - \bm{x}^{\star}\|^2 \right]
		-   \mathbb{E} \left[\|\bm{x}^{k+1} - \bm{x}^{\star} \|^2 \right]   + \mathcal{O}( K^{-1}).
	\end{split}
\end{equation*}
Summing the above inequality over $k$, we have
\begin{align*}
	\begin{split}
		&\,  \frac{\rho}{8m\kappa^2} \sum_{k=0}^{K-1} \mathbb{E} \left[ \|\bm{x}^k - \Pi_{C}(\bm{x}^k) \|^2 \right] \\
		\leq & \, \| \bm{x}^0 - \bm{x}^{\star} \|^2 - \mathbb{E} \left[ \|\bm{x}^{K} - \bm{x}^{\star}\|^2 \right] + \mathcal{O}(1) \sum_{k=0}^{K-1} \frac{1}{K} 
		\leq \mathcal{O}(1).
	\end{split}
\end{align*}
By the convexity of $\dist^2(\cdot,C)$, it follows that
\begin{equation*}
	\begin{split}
		&\,   \mathbb{E} \left[\dist^2(\bar{\bm{x}}^{K}, C) \right] \\
		\leq \, &\frac{1}{K} \sum_{k=0}^{K-1} \mathbb{E} \left[\|\bm{x}^k - \Pi_{C}(\bm{x}^k)\|^2 \right] 
		\leq  \mathcal{O} \left( \frac{1}{K} \right).
	\end{split}
\end{equation*}
Next, we prove the convergence rate of the optimality gap. By the definition of $\bm{v}^{k}$, we have
\begin{equation*}
	\begin{split}
		&\mathbb{E} \left[  \langle \bm{v}^{k}, \bm{x}^k - \bm{x}^{\star} \rangle \mid \mathcal{F}_{k} \right] 
		=\langle \nabla f(\bm{x}^k),\,\bm{x}^k-\bm{x}^{\star}\rangle  \\
		\geq&  f(\bm{x}^k) - f(\bm{x}^{\star}) = f(\bm{x}^k) - f^{\star} ,
	\end{split}
\end{equation*}
which, together with that $\bm{x}^\star\in H_k, C_0$ and Lemma~\ref{lema: SVRG-type-Property}, implies
\begin{equation*}
	\begin{split}
		&\, \mathbb{E} \left[ \|\bm{x}^{k+1} - \bm{x}^{\star} \|^2 \mid \mathcal{F}_k \right]\le \mathbb{E} \left[ \|\bm{x}^k - \alpha \bm{v}^k - \bm{x}^{\star} \|^2 \mid \mathcal{F}_k \right]\\
	= &\, \mathbb{E} \left[ \|\bm{x}^k - \bm{x}^\star\|^2 + \alpha^2 \|\bm{v}^k\|^2 - 2\alpha \langle \bm{x}^k - \bm{x}^\star , \bm{v}^k \rangle\mid \mathcal{F}_k \right]\\
		\leq &\,  ( 1 + 8L^2 \alpha^2 )\|\bm{x}^k - \bm{x}^{\star} \|^2 + 16L^2 \alpha^2 \| \tilde{\bm{x}}^{l} - \bm{x}^{\star} \|^2 \\
		&\, - 2 \alpha \left( f(\bm{x}^k) - f^{\star} \right)  + 4 \| \nabla f(\bm{x}^{\star}) \|^2 \alpha^2   + 12 R^2 \alpha^2  \\
		\le &\, \|\bm{x}^k - \bm{x}^{\star} \|^2 - 2  \alpha \left( f(\bm{x}^k) - f^{\star} \right) + \mathcal{O} \left( \frac{1}{K} \right),
	\end{split}
\end{equation*}
where the last inequality follows from the boundedness of $C_0$ and the fact that $\tilde{\bm{x}}^{l}, \bm{x}^{k}, \bm{x}^{\star}\in C_0$ and $\alpha = \mathcal{O}(1/\sqrt{K})$.
Taking expectation on both sides, we obtain
\begin{equation*}
	\begin{split}
		& \, \mathbb{E} \left[ f(\bm{x}^k) - f^{\star} \right] \\ 
		\leq &\, \frac{1}{2\alpha} \left(\mathbb{E} \left[ \|\bm{x}^{k} - \bm{x}^{\star} \|^2 \right] - \mathbb{E} \left[ \|\bm{x}^{k+1} - \bm{x}^{\star} \|^2 \right] \right) + \mathcal{O}\left( \frac{1}{\sqrt{K}}\right).
	\end{split}
\end{equation*}
Summing of the last inequality over $k$, we get
\begin{equation*}
	\begin{split}
		&\, \sum_{k=0}^{K-1}\mathbb{E} \left[ f(\bm{x}^k) - f^{\star} \right] \\
		\leq &\, \frac{1}{2\alpha}  \sum_{k=0}^{K-1} \left(\mathbb{E} \left[ \| \bm{x}^{k} - \bm{x}^{\star} \|^2 \right] - \mathbb{E} \left[ \| \bm{x}^{k+1} - \bm{x}^{\star} \|^2 \right] \right) \\
		&\, +\mathcal{O}(1)\sum_{k=0}^{K-1}\frac{1}{\sqrt{K}}  \\
		\leq &\, \mathcal{O}(1)  \sqrt{K+1}  +  \mathcal{O}(1)  \sqrt{K}
		=\mathcal{O}(\sqrt{K}),
	\end{split}
\end{equation*}
which, together with the convexity of $f$, yields
\begin{equation*}
	\begin{split}
		\mathbb{E} [f(\bar{\bm{x}}^{K}) - f^{\star}] \leq &\,  \frac{1}{K} \sum_{k=0}^{K-1} \mathbb{E}[f(\bm{x}^k) - f^{\star}] 
		\leq \, \mathcal{O} \left( \frac{1}{\sqrt{K}} \right).
	\end{split}
\end{equation*}
This completes the proof.

\subsection{Proof of Theorem~\ref{thm:var_redu_rate}}\label{subsec:ver_redu_rate}
To prove the feasibility gap, by Lemma~\ref{fea-gap-C-bounded}, we get
\begin{equation*} 
 \begin{split}
       &\frac{\rho}{8m\kappa^2}\mathbb{E}[ \dist^2(\bm{x}^k, C)] = \frac{\rho}{8m\kappa^2}\mathbb{E}[\|\bm{x}^k-\Pi_{C}(\bm{x}^k)\|^2] \\
\le\;&\mathbb{E}[\|\bm{x}^k-\Pi_{C}(\bm{x}^k)\|^2]-
\mathbb{E}[\|\bm{x}^{k+1}-\Pi_C(\bm{x}^{k+1})\|^2]
+\mathcal{O}(\alpha_k^2). 
\end{split}
\end{equation*}
Summing this inequality over $k$, we get 
\begin{equation}\label{ineq:proof:11}
\begin{split}
        &\,   \frac{1}{K}\sum_{k=0}^K\mathbb{E}[\|\bm{x}^k-\Pi_C(\bm{x}^k)\|^2] \\
\le &\, \mathcal{O}\left( \tfrac{1}{K} \right) + \mathcal{O}(\tfrac{1}{K} )  \sum_{k=1}^{K}\frac{1}{k^{\frac{2}{3}}} \le   \mathcal{O}\left(\tfrac{1}{K^{\frac{2}{3}}}\right).
    \end{split}
\end{equation}
The desired bound on feasibility gap then follows from the convexity of $\dist^2(\cdot, C)$.
For the optimality gap, since $R=0$, inequality~\eqref{ineq: Lip-cons-f} implies
\begin{equation*}
\begin{split}
f(\bm{x}^{k+1}) \le&\, f(\bm{x}^k)+\langle \nabla f(\bm{x}^k),\bm{x}^{k+1}-\bm{x}^k\rangle \\
&\,+\tfrac{L}{2}\|\bm{x}^{k+1}-\bm{x}^k\|^2.
\end{split}
\end{equation*}
Let $\bm{x}^\star\in X^\star$. By the convexity of $f$, we get
\begin{align*}
-f(\bm{x}^{\star})
\le\; -(f(\bm{x}^k)+\langle \nabla f(\bm{x}^k), \bm{x}^{\star}-\bm{x}^k \rangle ).
\end{align*}
Adding the last two inequalities and taking conditional expectation on the resulting inequality, we get
\begin{equation}\label{f_x+1-f_x^star-C-bounded}
\begin{split}        
&\, \mathbb{E}[f(\bm{x}^{k+1})-f(\bm{x}^{\star})\mid \mathcal{F}_k]\\
\le & \, \mathbb{E}[ \langle \bm{v}^k,  \bm{x}^{k+1}-\bm{x}^{\star}  \rangle
+\tfrac{L}{2}\|\bm{x}^{k+1}-\bm{x}^k\|^2 \mid \mathcal{F}_k ]\\
&\, +\mathbb{E}[\langle \nabla f(\bm{x}^k)- \bm{v}^k,  \bm{x}^{k+1}-\bm{x}^{\star} \mid \mathcal{F}_k \rangle] \\
= &\, \mathbb{E}[ \langle \bm{v}^k,  \bm{x}^{k+1}-\bm{x}^{\star}  \rangle
+\tfrac{L}{2}\|\bm{x}^{k+1}-\bm{x}^k\|^2 \mid \mathcal{F}_k ] \\
&\, +\mathbb{E}[\langle \nabla f(\bm{x}^k)- \bm{v}^k,  \bm{x}^{k+1}-\bm{x}^{k}  \rangle \mid \mathcal{F}_k ] \\
\le &\, \mathbb{E}[ \langle \bm{v}^k,  \bm{x}^{k+1}-\bm{x}^{\star}  \rangle + \tfrac{L}{2}\|\bm{x}^{k+1}-\bm{x}^k\|^2\mid \mathcal{F}_k ] \\
&\,+ \mathbb{E}[ 2L\| \bm{x}^{k+1}-\bm{x}^{k} \|^2   +  \tfrac{2}{L}\|\nabla f(\bm{x}^k)- \bm{v}^k \|^2\mid \mathcal{F}_k ],
    \end{split}
\end{equation}
where the equality follows from $\mathbb{E}[\langle \nabla f(\bm{x}^k)- \bm{v}^k,  \bm{x}^{k}-\bm{x}^{\star}  \rangle
\mid \mathcal{F}_k]=0$.
By the definition of $\bm{v}^k$, we have
\begin{equation}
\label{vk-nabla(f(xk))}
\begin{split}
&\, \mathbb{E}[ \|\nabla f(\bm{x}^k)- \bm{v}^k \|^2 \mid \mathcal{F}_k]\\
\le & \, 2\mathbb{E}[ \|\tfrac{1}{b}\sum_{i\in I_k}
(\nabla f_{i}(\bm{x}^{k})-\nabla f_{i}(\tilde{\bm{x}}^{l})) \|^2\mid \mathcal{F}_k] \\
&\, +2 \mathbb{E}[\|\nabla f(\tilde{\bm{x}}^{l}) -\nabla f(\bm{x}^k)\|^2\mid \mathcal{F}_k]\\
\le  &\, 4L\mathbb{E}[\|\bm{x}^{k}-\tilde{\bm{x}}^{l}\|^2\mid \mathcal{F}_k].
\end{split}
\end{equation}
Substituting~\eqref{vk-nabla(f(xk))} into \eqref{f_x+1-f_x^star-C-bounded} yields
\begin{equation}\label{ineq:fxk-fxsol}
\begin{split}
&\,\mathbb{E}[f(\bm{x}^{k+1})-f(\bm{x}^{\star})\mid \mathcal{F}_k]\\
\le&\, \mathbb{E}[ \langle \bm{v}^k,  \bm{x}^{k+1}-\bm{x}^{\star}  \rangle + \tfrac{5L}{2}\|\bm{x}^{k+1}-\bm{x}^k\|^2 \mid \mathcal{F}_k]\\
&\, + 8\mathbb{E}[\|\bm{x}^{k}-\tilde{\bm{x}}^{l}\|^2\mid \mathcal{F}_k].
    \end{split}
\end{equation}
We then claim that
\begin{equation}\label{ineq:vk-xk-xsol}
\begin{split}
& \, \mathbb{E}[ \langle \bm{v}^k,  \bm{x}^{k+1}-\bm{x}^{\star}  \rangle + \tfrac{5L}{2}\|\bm{x}^{k+1}-\bm{x}^k\|^2\mid \mathcal{F}_k] \\
\le &\, \tfrac{1}{2\alpha_k}\mathbb{E}[(\|\bm{x}^k-\bm{x}^{\star}\|^2-\|\bm{x}^{k+1}-\bm{x}^{\star}\|^2)\mid \mathcal{F}_k].
\end{split}
\end{equation}
Indeed, if $\bm{\xi}^k=\bm{0}$ or $(\phi_{j_k}(\bm{x}^k)-\alpha_k\langle \bm{\xi}^k, \bm{v}^{k}\rangle )_{+}=0$, from the definition of $\bm{x}^{k+1}$, we have $\bm{v}^k=\tfrac{1}{\alpha_k}(\bm{x}^{k}-\bm{x}^{k+1})$
and
\begin{equation*}
\begin{split}
&\, \mathbb{E}[ \langle \bm{v}^k,  \bm{x}^{k+1}-\bm{x}^{\star}  \rangle + \tfrac{5L}{2}\|\bm{x}^{k+1}-\bm{x}^k\|^2\mid \mathcal{F}_k]\\
= &\,\mathbb{E}[\tfrac{1}{2\alpha_k} (\|\bm{x}^k-\bm{x}^{\star}\|^2-\|\bm{x}^{k+1}-\bm{x}^{\star}\|^2-\|\bm{x}^{k+1}-\bm{x}^k\|^2) \\
&\,+ \tfrac{5L}{2}\|\bm{x}^{k+1}-\bm{x}^k\|^2\mid \mathcal{F}_k]\\
\le &\, \tfrac{1}{2\alpha_k}\mathbb{E}[ \|\bm{x}^k-\bm{x}^{\star}\|^2-\|\bm{x}^{k+1}-\bm{x}^{\star}\|^2 \mid \mathcal{F}_k],
\end{split}
\end{equation*}
where the inequality follows from $5L/2\le 1/(2\alpha_k)$.
If $\bm{\xi}^k\neq\bm{0}$ and
$(\phi_{j_k}(\bm{x}^k)-\alpha_k\langle \bm{\xi}^k, \bm{v}^{k}\rangle )_{+} > 0$,
we have 
\[
\bm{v}^k=\tfrac{1}{\alpha_k} \left(\bm{x}^{k}-\bm{x}^{k+1}-\tfrac{ \left( \phi_{j_k}(\bm{x}^k)-\alpha_k\langle \bm{\xi}^k, \bm{v}^{k}\rangle \right) }{\|\bm{\xi}^k\|^2}\bm{\xi}^k \right) 
\]
and 
\begin{equation}\label{ineq:vk-xk-x}
\begin{split}
&\, \langle \bm{v}^k,  \bm{x}^{k+1}-\bm{x}^{\star}  \rangle
+\tfrac{5L}{2}\|\bm{x}^{k+1}-\bm{x}^k\|^2\\
\le &\, \tfrac{1}{2\alpha_k}(\|\bm{x}^k-\bm{x}^{\star}\|^2-\|\bm{x}^{k+1}-\bm{x}^{\star}\|^2)\\
&\, -  \tfrac{ \left( \phi_{j_k}(\bm{x}^k)-\alpha_k\langle \bm{\xi}^k, \bm{v}^{k}\rangle \right)  }{\alpha_k \|\bm{\xi}^k\|^2}\langle \bm{\xi}^k,
\bm{x}^{k+1}-\bm{x}^{\star}\rangle.
\end{split}
\end{equation}
By the convexity of $\phi_{j_k}$, we have  
\begin{equation*}
\begin{split}
&\, \langle -\bm{\xi}^k, \bm{x}^{k+1}-\bm{x}^{\star}\rangle
=  \langle \bm{\xi}^k, \bm{x}^{\star}-\bm{x}^{k}+\bm{x}^{k}-\bm{x}^{k+1}\rangle\\
\le &\, \phi_{j_k}(\bm{x}^{\star})-\phi_{j_k}(\bm{x}^{k})
+ \langle \bm{\xi}^k, \bm{x}^{k}-\bm{x}^{k+1}\rangle = \phi_{j_k}(\bm{x}^{\star}) \le 0,
    \end{split}
\end{equation*}
which, together with \eqref{ineq:vk-xk-x} and $\tfrac{ \left( \phi_{j_k}(\bm{x}^k)-\alpha_k\langle \bm{\xi}^k, \bm{v}^{k}\rangle \right)  }{\alpha_k \|\bm{\xi}^k\|^2}>0$, implies \eqref{ineq:vk-xk-xsol}.
Substituting  (\ref{ineq:vk-xk-xsol}) into (\ref{ineq:fxk-fxsol}), we obtain
\begin{equation*}
\begin{split}
& \, \mathbb{E}[f(\bm{x}^{k+1})-f(\bm{x}^{\star})\mid \mathcal{F}_k  ]\\
\le &\, \tfrac{1}{2\alpha_k}\mathbb{E}[(\|\bm{x}^k-\bm{x}^{\star}\|^2-
\|\bm{x}^{k+1}-\bm{x}^{\star}\|^2) \mid \mathcal{F}_k]
+ 8\mathbb{E}[\| \bm{x}^{k}-\tilde{\bm{x}}^l\|^2\mid \mathcal{F}_k],
\end{split}
\end{equation*}
which, upon taking expectation and using Lemma~\ref{x^k+1-tilde_x^l}, yields
\begin{equation*}
\begin{split}
&\, \mathbb{E}[f(\bm{x}^{k+1})-f(\bm{x}^{\star})]\\
\le &\, \tfrac{1}{2}(\tfrac{1}{\alpha_k}\mathbb{E}[\|\bm{x}^k-\bm{x}^{\star}\|^2]
-\tfrac{1}{\alpha_{k+1}}\mathbb{E}[\|\bm{x}^{k+1}-\bm{x}^{\star}\|^2]) \\
&\, +\tfrac{1}{2}(\tfrac{1}{\alpha_{k+1}}-\tfrac{1}{\alpha_{k}})
\mathbb{E}[\|\bm{x}^{k+1}-\bm{x}^{\star}\|^2]\\
&\, +\mathcal{O}(\alpha_k^2)
+16r \sum_{t=lr}^{t=k} \mathbb{E}[\|\bm{x}^t-\Pi_C(\bm{x}^t)\|^2].
\end{split}
\end{equation*}
Summing this inequality over $k$ and using Lemma~\ref{lem:x^k_bound},  we have
\begin{equation*}
\begin{split}
&\, \sum_{k=0}^{K-1} (\mathbb{E}[f(\bm{x}^{k+1})]-f(\bm{x}^{\star}))\\
\le &\,\tfrac{1}{2\alpha_0}\|\bm{x}^0-\bm{x}^{\star}\|^2
+\tfrac{1}{\alpha_{K+1}}\mathcal{O}(1)+\sum_{k=0}^{K - 1}\mathcal{O}(\alpha_k^2) \\
&\, +16 r^2 \sum_{k=0}^{K-1} \mathbb{E}[\|\bm{x}^k-\Pi_C(\bm{x}^k)\|^2].
\end{split}
\end{equation*}
The last inequality, the definition of $\bar{\bm{x}}$ and~\eqref{ineq:proof:11} then imply 
\begin{equation*}
\begin{split}
&\, \mathbb{E}[f(\bar{\bm{x}}^K)]-f(\bm{x}^{\star})\\
\le &\, (1 +\tfrac{1}{\alpha_{K+1}})\mathcal{O}(\tfrac{1}{K}) + \frac{1}{K}\sum_{k=0}^K\mathcal{O}(\alpha_k^2) + \mathcal{O}\left(\tfrac{1}{K^{\frac{2}{3}}}\right) \\
\le &\, \mathcal{O}(\tfrac{1}{K}) +\mathcal{O}(\tfrac{K^{\frac{1}{3}}}{K}) +\frac{1}{K}\sum_{k=0}^K\mathcal{O}(\alpha_k^2)+\mathcal{O}\left(\tfrac{1}{K^{\frac{2}{3}}}\right) \\
\le &\, \mathcal{O}(\tfrac{1}{K}) + \mathcal{O}(\tfrac{1}{K^{\frac{2}{3}}}) \le 
\mathcal{O}\left( \tfrac{1}{K^{\frac{2}{3}}} \right).
\end{split}
\end{equation*}
This completes the proof.

\subsection{Proof of Theorem~\ref{rate-qg}}\label{subsec:converge_rate_QG}
By Theorem~\ref{convergence1}, $\{\bm{x}^{k}\}$ converges almost surely to a point in $X^{\star}$.
Therefore, the sequence $\{ \mathbb{E}[\|\bm{x}^{k}\| ]\}_k$ is bounded, and by Lemma~\ref{lema: Proj-Property}\ref{lema: proj-property-nonexpan}, inequality~\eqref{ineq: Lip-cons-f} and triangle inequality, so are the sequences $\{\mathbb{E}[ \|\tilde{\bm{x}}^{l}\|] \}_l$, $\{\mathbb{E}[ \|\Pi_{X^{\star}} (\bm{x}^{k})\|]\}_k$, $\{ \mathbb{E}[\|\bm{x}^{k} - \Pi_{X^{\star}} (\bm{x}^{k}) \|]\}_k $, $\{\mathbb{E}[\|\bm{x}^{k} - \Pi_C(\bm{x}^{k})\| ]\}_k$ and $\{\mathbb{E} \left[  \|\nabla f (\Pi_{X^{\star}} (\bm{x}^k))\|^2 \right] \}_k$.
Using Lemma~\ref{lema: Proj-Property}\ref{lema: proj-property-nonexpan}, 
\begin{equation*}
	\begin{split}
		& \,   \|\bm{x}^{k} -   \Pi_{X^{\star}} (\bm{x}^{k}) \|  \\
		\leq &\,  \| \bm{x}^{k}  - \Pi_C (\bm{x}^{k} ) \| + \| \Pi_C (\bm{x}^{k} ) - \Pi_{X^{\star}} (\Pi_C (\bm{x}^{k} ) ) \| \\
		&\, + \| \Pi_{X^{\star}} (\Pi_C (\bm{x}^{k} ) ) - \Pi_{X^{\star}} (\bm{x}^{k}) \| \\
		\leq  & \, 2\| \bm{x}^{k}  - \Pi_C (\bm{x}^{k} ) \| +  \| \Pi_C (\bm{x}^{k} ) - \Pi_{X^{\star}} (\Pi_C (\bm{x}^{k} ) ) \|,
	\end{split}
\end{equation*} 
which implies that 
\begin{equation*}
	\begin{split}
		&\, -\| \Pi_C (\bm{x}^{k} ) - \Pi_{X^{\star}} (\Pi_C (\bm{x}^{k} ) ) \|^2 \\
		\leq     &\, 4\| \bm{x}^{k}  - \Pi_C (\bm{x}^{k} ) \|^2 -  \frac{1}{2}\|\bm{x}^{k} -   \Pi_{X^{\star}} (\bm{x}^{k}) \|^2.
	\end{split}
\end{equation*}
Using this inequality, Assumption~\ref{ass:QG} and the fact that $\alpha_k\geq \frac{8}{\nu (k+1)}$, we have
\begin{align*}
	&\, -\alpha_k (f(\Pi_C (\bm{x}^k)) - f^{\star}) \\
	\leq   &\,  -\frac{\alpha_k\nu}{2}\|\Pi_C(\bm{x}^k) - \Pi_{X^{\star}}\left( \Pi_C(\bm{x}^k) \right)    \|^2  \\ 
	\leq  & \,  2\nu \alpha_k   \|\bm{x}^k - \Pi_{C}(\bm{x}^k)\| ^2-\frac{2}{k+1} \|\bm{x}^{k} -  \Pi_{X^{\star}} (\bm{x}^{k})\|^2.
\end{align*}
This inequality and Lemma~\ref{lemma: one-iteration} with $\lambda=32m\kappa^2\rho^{-1}$ yield that
\begin{align} 
		&\, \mathbb{E}\left[\|\bm{x}^{k+1}-\Pi_{X^{\star}} (\bm{x}^{k+1})\|^2\mid\mathcal{F}_{k}\right] \notag  \\  
		\leq &\,  \mathbb{E}\left[\|\bm{x}^{k+1}- \Pi_{X^{\star}} (\bm{x}^{k})\|^2\mid\mathcal{F}_{k}\right]  \notag  \\
		\leq\, & (1 + \alpha_k^2(24 L^2+L^2\lambda))\|\bm{x}^k-  \Pi_{X^{\star}} (\bm{x}^{k})\|^2 \notag  \\
		&\, + 48 L^2 \alpha_k^2\|\tilde{\bm{x}}^{l}- \Pi_{X^{\star}} (\bm{x}^{k})\|^2 + 36 R^2\alpha_k^2  \notag  \\
		& + \alpha_k^2\left(2\lambda R^2+ (\lambda+12) \|\nabla f(  \Pi_{X^{\star}} (\bm{x}^{k}))\|^2\right)\notag  \\
		& -(\tfrac{3\rho}{8 m \kappa^2} - 2\alpha_k L ) \| \bm{x}^k-\Pi_{C}( \bm{x}^k ) \|^2 \notag \\
		&\, - \left(\alpha_k +\frac{8}{\nu (k+1)} \right) \left( f( \Pi_C(\bm{x}^k)) - f^{\star} \right)  \label{ineq:thm3.3-2-1-2}   \\
		\leq & \, \left(1 + \alpha_k^2(24 L^2+L^2\lambda) -\frac{2}{k+1} \right)\|\bm{x}^k-  \Pi_{X^{\star}} (\bm{x}^{k}) \|^2 \notag  \\
		& \, -  \frac{8}{\nu (k+1)}\left( f( \Pi_C(\bm{x}^k)) - f^{\star} \right) \notag  \\
		&\,
		- \left( \frac{3\rho}{8m\kappa^2} - 2\alpha_k  (L+\nu) \right)\| \bm{x}^k-\Pi_{C}( \bm{x}^k ) \|^2 \notag  \\
		& \, + \alpha_k^2\left(2\lambda R^2+ (\lambda+12) \|\nabla f(\Pi_{X^{\star}} (\bm{x}^{k}))\|^2 + 36R^2 \right) \notag  \\ 
		&\, + 48 L^2 \alpha_k^2\|\tilde{\bm{x}}^{l}- \Pi_{X^{\star}} (\bm{x}^{k})\|^2.\notag 
\end{align} 
Noting that
\begin{equation*}
	\begin{split}
		&\,    \frac{3\rho}{8m\kappa^2} - 2\alpha_k (L+\nu)    \\
		= &\,    \frac{3\rho}{8m\kappa^2} - 2\alpha_k (L+\nu)  + 4\alpha_k^2 (L+\nu)^2 m \kappa^2 \rho^{-1}  - \mathcal{O}(\alpha_k^2)  \\
		= &\, \frac{\rho}{8m\kappa^2} +\left( \sqrt{\frac{\rho}{4m\kappa^2}} - \sqrt{4m\kappa^2\rho^{-1}} \alpha_k (L+\nu)\right)^2 - \mathcal{O}(\alpha_k^2)\\
		\geq &\, \frac{\rho}{8m\kappa^2}- \mathcal{O}(\alpha_k^2),
	\end{split}
\end{equation*}
and taking total expectation on \eqref{ineq:thm3.3-2-1-2}, we have 
\begin{equation*} \label{ineq:thm3.3-3}
	\begin{split} 
		& \mathbb{E}\left[\|\bm{x}^{k+1}- \Pi_{X^{\star}} (\bm{x}^{k+1})\|^2 \right]  \\
		\leq &  \left(1 - \frac{2}{k+1} \right)\mathbb{E}[ \|\bm{x}^k- \Pi_{X^{\star}} (\bm{x}^{k})\|^2] \\
		&
		-   \frac{8}{\nu(k+1)}\mathbb{E}\left[  f( \Pi_C(\bm{x}^k)) -  f^{\star} \right]  +  \mathcal{O}( \alpha_k^2)
		\\
		&   - \left( \frac{3\rho}{8m\kappa^2} - 2\alpha_k (L+\nu)  \right)\mathbb{E}[\| \bm{x}^k-\Pi_{C}( \bm{x}^k ) \|^2] \\
		\leq &   \left(1 - \frac{2}{k+1}\right)\mathbb{E} \left[\|\bm{x}^k- \Pi_{X^{\star}} (\bm{x}^{k})\|^2\right] \\
		&
		-  \frac{8}{\nu(k+1)} \mathbb{E}\left[  f( \Pi_C(\bm{x}^k)) - f^{\star} \right] 
		\\
		& -  \frac{\rho}{8m \kappa^2} \mathbb{E}\left[ \| \bm{x}^k-\Pi_{C}( \bm{x}^k ) \|^2\right]+  \mathcal{O}(\alpha_k^2),
	\end{split}
\end{equation*}
where we used the boundedness of the sequences mentioned at the beginning of the proof.
Multiplying both sides of the last inequality by~$(k+1)k$ and using the fact that $\alpha_k = \mathcal{O}(1/k)$, we have that for all $K\geq k$,
\begin{equation*} 
	\begin{split}
		&\, \frac{\rho}{8m\kappa^2} (k+1)k\mathbb{E}\left[\| \bm{x}^k-\Pi_{C}( \bm{x}^k ) \|^2 \right] \\
		&\, + \frac{8}{\nu  (K+1)} (k+1) k\mathbb{E}\left[ f( \Pi_C(\bm{x}^k)) - f^{\star} \right] \\
		\leq  &\,  (k-1)k\mathbb{E}\left[ \|\bm{x}^k- \Pi_{X^{\star}} (\bm{x}^{k})\|^2 \right] \\
		& \, -  (k+1)k\mathbb{E}\left[\|\bm{x}^{k+1}- \Pi_{X^{\star}} (\bm{x}^{k+1})\|^2\right] + \mathcal{O}(1).
	\end{split}
\end{equation*}
Summing the above inequality over $k$, we get
\begin{equation*} 
	\begin{split}
		&\,  \sum_{k=1}^{K}(k+1)k\left( \frac{\rho}{8m\kappa^2} \mathbb{E}\left[\|\bm{x}^k-\Pi_{C}( \bm{x}^k)\|^2\right]\right) \\
		+ &\,\sum_{k=1}^{K}(k+1)k\left( \frac{8}{\nu  (K+1)}   \mathbb{E}\left[  f( \Pi_C(\bm{x}^k)) - f^{\star} \right] \right) \\
		\leq &\,    \mathcal{O}\left(K\right),
	\end{split}
\end{equation*}
which, together with the convexity of $\| \cdot\|^2$ and $f$, implies that
\begin{equation}\label{ineq:thm3.3-5}
	\begin{split}
		&  \mathbb{E}\left[\|\bar{\bm{x}}^K - \hat{\bm{x}}^{K}\|^2\right] \leq \frac{K }{S_K}   \mathcal{O}(1) \,\,\,\, {\rm and}\\
		&
		\mathbb{E}\left[ f(\hat{\bm{x}}^{K}) - f^{\star}\right] \leq \frac{K (K+1)}{S_K} \mathcal{O}(1),
	\end{split}
\end{equation}
where $ \bar{\bm{x}}^{K} = \frac{1}{S_K} \sum_{k\in [K]}k (k+1)\bm{x}^{k}$, $ \hat{\bm{x}}^{K} = \frac{1}{S_K} \sum_{k\in [K]}k(k+1)\Pi_C(\bm{x}^{k})$ and
\begin{equation*}
	\begin{split}
		& S_K = \sum_{k\in [K]}(k^2+k) 
		\ge  \frac{1}{6} K^3.
	\end{split}
\end{equation*}
It then follows from~\eqref{ineq:thm3.3-5} and $S_K\ge K^3/6$ that  
\begin{equation*}\label{ineq:thm3.3-5-1}
	\begin{split}
		& \mathbb{E}\left[\dist^2(\bar{\bm{x}}^{K}, C)\right] \leq \, \mathbb{E}\left[\|\bar{\bm{x}}^{K} - \hat{\bm{x}}^{K}\|^2\right] \leq \mathcal{O}\left(\frac{1}{K^2}\right)\\
		&\text{and }\mathbb{E}\left[ f(\hat{\bm{x}}^{K}) - f^{\star}\right] \leq   \mathcal{O}\left(\frac{1}{K}\right).
	\end{split}
\end{equation*}
These two inequalities and the mean value theorem imply the existence of $\theta\in [0,1]$ such that
\begin{equation*}\label{ineq:thm3.3-6}
	\begin{split}
		&\, \mathbb{E}\left[ f(\bar{\bm{x}}^{K}) - f^{\star}  \right]\le	\mathbb{E}\left[ |f(\bar{\bm{x}}^{K}) - f^{\star} | \right] \\
		\leq   &\,   \mathbb{E} \left[ f(\hat{\bm{x}}^{K}) - f^{\star}\right]  +   \mathbb{E}\left[|f(\bar{\bm{x}}^{K})- f(\hat{\bm{x}}^{K})|\right] \\
		\leq &\, \mathcal{O}\left(\frac{1}{K}\right) + \mathbb{E}\left[ \|\nabla f(\bar{\bm{x}}^{K} + \theta (\hat{\bm{x}}^{K} - \bar{\bm{x}}^{K}))\| \|\bar{\bm{x}}^{K}-\hat{\bm{x}}^{K}\|\right]\\
		\leq &\,  \mathcal{O}\left(\frac{1}{K}\right) + \frac{K}{2}\mathbb{E}\left[  \|\bar{\bm{x}}^{K}-\hat{\bm{x}}^{K}\|^2\right]\\
		&\, +\frac{1}{2K} \mathbb{E}\left[  \|\nabla f(\bar{\bm{x}}^{K} + \theta (\hat{\bm{x}}^{K} - \bar{\bm{x}}^{K}))\|^2\right]\\
		\leq &\,  \mathcal{O}\left(\frac{1}{K}\right) +  \mathcal{O}\left(\frac{1}{K}\right) \mathbb{E}\left[ \|\nabla f(\bar{\bm{x}}^{K} + \theta (\hat{\bm{x}}^{K} - \bar{\bm{x}}^{K}))\|^2\right].
	\end{split}
\end{equation*}
The desired result follows from the boundedness of the sequence $\{ \mathbb{E} \left[ \|\nabla f(\bar{\bm{x}}^{K} + \theta (\hat{\bm{x}}^{K} - \bar{\bm{x}}^{K}))\|^2 \right]\}_K $, which can be proved as follows:
\begin{equation*}\label{ineq:thm3.3-7}
	\begin{split}
		&  \mathbb{E} \left[ \|\nabla f(\bar{\bm{x}}^{K} + \theta (\hat{\bm{x}}^{K} - \bar{\bm{x}}^{K}))\|^2 \right]\\
		\leq  &\,  2 \mathbb{E} \left[  \|\nabla f(\bar{\bm{x}}^{K} + \theta (\hat{\bm{x}}^{K} - \bar{\bm{x}}^{K}))-\nabla f (\bm{x}^{0}) \|^2 \right]  +2\|\nabla f (\bm{x}^{0}) \|^2 \\
		\leq &\, 2L^2\mathbb{E}\left[ (\|\bar{\bm{x}}^{K} + \theta (\hat{\bm{x}}^{K} - \bar{\bm{x}}^{K}) - \bm{x}^{0}\| +1)^2\right]  +2\|\nabla f (\bm{x}^{0}) \|^2\\
		\leq &\, \mathcal{O}(1) + 8L^2 \mathbb{E}\left[ (\|\bar{\bm{x}}^{K} -  \bm{x}^{0}\|^2 + \|\hat{\bm{x}}^{K} - \bar{\bm{x}}^{K}\|^2) \right] \\
		\leq &\, \mathcal{O}(1) +16L^2 \|\bm{x}^{0}\|^2   + \mathcal{O}\left(\frac{1}{K}\right)\\
		&\,  +\frac{ 16L^2 }{S_K} \sum_{k\in [K]} (k+1)k \mathbb{E}\left[\| \bm{x}^{k}\|^2\right]
		\leq  \, \mathcal{O}(1),
	\end{split}
\end{equation*}
where the second inequality follows from~(\ref{ineq: Lip-cons-f}), the third from $\theta\in [0,1]$, the fourth from the convexity of $\|\cdot\|^2$, and the last from  the boundedness of the sequence $\left\{ \mathbb{E}\left[\| \bm{x}^{k}\|^2 \right] \right\}_k$.
This completes the proof.

\bibliographystyle{ieeetr}
\bibliography{refs}

\end{document}